\documentclass[12pt,reqno]{amsart} 
\usepackage{amssymb}
\usepackage{mathrsfs}
\usepackage{enumerate}
\usepackage[all]{xy}
\newcommand{\mynote}[1]{\noindent\textbf{[#1]}}

\newcommand{\4}[1]{\widebar{#1}}
\newcommand{\9}[1]{{}^{#1}\!}

\let\oldcirc=\circ
\renewcommand{\circ}{\mathchoice
    {\mathbin{\scriptstyle\oldcirc}}{\mathbin{\scriptstyle\oldcirc}}
    {\mathbin{\scriptscriptstyle\oldcirc}}
    {\mathbin{\scriptscriptstyle\oldcirc}}}

\DeclareMathAlphabet\EuR{U}{eur}{m}{n}
\SetMathAlphabet\EuR{bold}{U}{eur}{b}{n}

\def\II(#1,#2){\cals(#2)_{\ge#1}}

\SelectTips{cm}{10} \UseTips   
\newlength{\upto}\newlength{\dnto}
\numberwithin{equation}{section}

\mathcode`\:="003A
\mathchardef\cdot="0201

\newcommand{\newsect}[1]{\bigskip\section{#1}\setcounter{table}{0}}

\renewenvironment{enumerate}[1][]
{\begin{enumerat}[#1]\setlength{\itemsep}{6pt}}{\end{enumerat}}

\newenvironment{enuma}{\begin{enumerate}[{\rm(a) }]}{\end{enumerate}}
\newenvironment{enum1}[1][]{\begin{enumerate}[{\rm({#1}1) }]}{\end{enumerate}}

\renewenvironment{itemize}
{\begin{itemiz}\setlength{\itemsep}{6pt}\setlength{\itemindent}{-20pt}}
{\end{itemiz}}

\def\beq#1\eeq{\begin{equation*}#1\end{equation*}}
\def\beqq#1\eeqq{\begin{equation}#1\end{equation}}

\let\emptyset=\varnothing

\newcommand{\longline}{\bigskip\centerline{\hbox to 5cm{\hrulefill}}\bigskip}

\newcommand{\widebar}[1]{\overset{\mskip3mu\hrulefill\mskip3mu}{#1}
		\vphantom{#1}}

\renewcommand{\:}{\colon}

\newcommand{\orb}{\mathcal{O}}

\newcounter{let} 
\loop\stepcounter{let}
\expandafter\edef\csname cal\alph{let}\endcsname%
{\noexpand\mathcal{\Alph{let}}}
\ifnum\thelet<26\repeat

\loop\stepcounter{let}
\expandafter\edef\csname\Alph{let}\alph{let}\endcsname%
{\noexpand\mathcal{\Alph{let}}}
\ifnum\thelet<26\repeat

\renewcommand{\cals}{\mathscr{S}}

\newcommand{\tdef}[2][]{\expandafter\newcommand\csname#2\endcsname%
{#1\textup{#2}}}
\tdef{Iso}   \tdef{Out}    \tdef{Aut}    \tdef{Inn}    \tdef{Hom}
\tdef{Inj}   \tdef{map}    \tdef{Ker}    \tdef{Ob}     \tdef{End}
\tdef{Mor}   \tdef{Res}    \tdef{Spin}   \tdef{Id}     \tdef{Fr}
\tdef{rk}    \tdef{conj}   \tdef{incl}   \tdef{proj}   \tdef{diag} 
\tdef{trf}   \tdef{Sol}    \tdef{He}     \tdef{Sz}     \tdef{cj}
\tdef{pr}    \tdef{Baum}   \tdef{Syl}    \tdef{Comp}   \tdef{Rep}
\tdef{GL}    \tdef{SL}     \tdef{Stab}   \tdef{Br}     \tdef{PSL} 
\tdef{Mat}   \tdef{id}     \tdef{br}     \tdef{Fix}    \tdef{IPr}
\tdef{Irr}   \tdef{res}    \tdef{Ext}
\tdef{supp}  \tdef[_]{typ}   \tdef[^]{op}    \tdef[^]{ab}

\newcommand{\fdef}[1]{\expandafter\newcommand\csname#1\endcsname%
{\mathfrak{#1}}}
\fdef{X}  \fdef{red}  \fdef{foc}  \fdef{hyp} \fdef{N}

\newcommand{\bbdef}[1]{\expandafter\newcommand%
\csname#1\endcsname{\mathbb{#1}}}
\bbdef{C} \bbdef{F} \bbdef{G} \bbdef{R} \bbdef{Z}

\newcommand{\gen}[1]{\langle{#1}\rangle}
\newcommand{\Gen}[1]{\bigl\langle{#1}\bigr\rangle}

\let\nsg=\normal

\newcommand{\syl}[2]{\textup{Syl}_{#1}(#2)}
\newcommand{\sylp}[1]{\syl{p}{#1}}
\newcommand{\autf}{\Aut_{\calf}}

\newcommand{\outf}{\Out_{\calf}}

\newcommand{\homf}{\Hom_{\calf}}

\newcommand{\sminus}{\smallsetminus}
\newcommand{\defeq}{\overset{\textup{def}}{=}}

\newcommand{\zploc}{\Z_{(p)}}

\newcommand{\pcom}{{}^\wedge_p}

\let\til=\widetilde

\newtheorem{Thm}{Theorem}[section]
\newtheorem{Prop}[Thm]{Proposition}
\newtheorem{Cor}[Thm]{Corollary}
\newtheorem{Lem}[Thm]{Lemma}

\newtheorem{Th}{Theorem}

\newcommand{\longleft}[1]{\;{\leftarrow%
\count255=0 \loop \mathrel{\mkern-6mu}%
    \relbar\advance\count255 by1\ifnum\count255<#1\repeat}\;}
\newcommand{\longright}[1]{\;{\count255=0 \loop \relbar\mathrel{\mkern-6mu}%
    \advance\count255 by1\ifnum\count255<#1\repeat\rightarrow}\;}
\newcommand{\Right}[2]{\overset{#2}{\longright#1}}
\newcommand{\RIGHT}[3]{\mathrel{\mathop{\kern0pt\longright#1}
	\limits^{#2}_{#3}}}

\newcommand{\LEFT}[3]{\mathrel{\mathop{\kern0pt\longleft#1}\limits^{#2}_{#3}}
}
\newcommand{\longleftright}[1]{\;{\leftarrow\mathrel{\mkern-6mu}%
    \count255=0\loop\relbar\mathrel{\mkern-6mu}%
    \advance\count255 by1\ifnum\count255<#1\repeat\rightarrow}\;} 
 
\newcommand{\onto}[1]{\;{\count255=0 \loop \relbar\joinrel
    \advance\count255 by1
    \ifnum\count255<#1 \repeat \twoheadrightarrow}\;}

\newcommand{\RLEFT}[3]{\mathrel{%
   \mathop{\vcenter{\baselineskip=0pt\hbox{$\kern0pt\longright#1$}%
   \hbox{$\kern0pt\longleft#1$}}}\limits^{#2}_{#3}}}
\newcommand{\RRIGHT}[3]{\mathrel{%
   \mathop{\vcenter{\baselineskip=0pt\hbox{$\kern0pt\longright#1$}%
   \hbox{$\kern0pt\longright#1$}}}\limits^{#2}_{#3}}}

\newcommand{\hclim}[1]{\setbox1=\hbox{\rm hocolim}
	\mathop{\vtop{\baselineskip=5pt\box1}}_{#1}}

\newcommand{\curs}{\EuR}
\newcommand{\Ab}{\curs{Ab}}

\newcommand{\spaces}{\curs{Top}}

\newcommand{\hotop}{\curs{hoTop}}
\renewcommand{\mod}{\mbox{-}\curs{mod}}

\theoremstyle{definition}

\newtheorem{Defi}[Thm]{Definition}

\theoremstyle{remark}

\newcommand{\higherlim}[2]{\displaystyle\setbox1=\hbox{\rm lim}
	\setbox2=\hbox to \wd1{\leftarrowfill} \ht2=0pt \dp2=-1pt
	\setbox3=\hbox{$\scriptstyle{#1}$}
	\def\test{#1}\ifx\test\empty
	\mathop{\mathop{\vtop{\baselineskip=5pt\box1\box2}}}\nolimits^{#2}
	\else
	\ifdim\wd1<\wd3
	\mathop{\hphantom{^{#2}}\vtop{\baselineskip=5pt\box1\box2}^{#2}}_{#1}
	\else
	\mathop{\mathop{\vtop{\baselineskip=5pt\box1\box2}}_{#1}}%
	\nolimits^{#2}
	\fi\fi}

\title{Existence and uniqueness of linking systems:\\ 
Chermak's proof via obstruction theory}

\author{Bob Oliver}

\address{LAGA, Institut Galil\'ee, Av. J-B Cl\'ement, 93430
Villetaneuse, France}
\email{bobol@math.univ-paris13.fr}
\thanks{B. Oliver was partially supported by the DNRF through a visiting 
professorship at the Centre for Symmetry and Deformation in Copenhagen; and 
also by UMR 7539 of the CNRS and by project ANR BLAN08-2\_338236, HGRT}


\newcommand{\LL}[3][*]{L^{#1}(#2;#3)}

\begin{document}

\begin{abstract} 
We present a version of a proof by Andy Chermak of the existence and 
uniqueness of centric linking systems associated to arbitrary saturated 
fusion systems.  This proof differs from the one in \cite{Chermak} in that 
it is based on the computation of higher derived functors of certain 
inverse limits.  This leads to a much shorter proof, but one which is 
aimed mostly at researchers familiar with homological algebra.
\end{abstract}

\maketitle

\bigskip

One of the central questions in the study of fusion systems is whether to 
each saturated fusion system one can associate a centric linking system, 
and if so, whether it is unique.  This question was recently answered 
positively by Andy Chermak \cite{Chermak}, using direct constructions.  His 
proof is quite lengthy, although some of the structures developed there 
seem likely to be of independent interest.

There is also a well established obstruction theory for studying this 
problem, involving higher derived functors of certain inverse limits.  This 
is analogous to the use of group cohomology as an ``obstruction theory'' 
for the existence and uniqueness of group extensions.  By using this 
theory, Chermak's proof can be greatly shortened, in part because it allows 
us to focus on the essential parts of Chermak's constructions, and in part 
by using results which are already established.  The purpose of this paper 
is to present this shorter version of Chermak's proof, a form which we hope 
will be more easily accessible to researchers with a background in topology 
or homological algebra.

A \emph{saturated fusion system} over a finite $p$-group $S$ is a category 
whose objects are the subgroups of $S$, and whose morphisms are certain 
monomorphisms between the subgroups.  This concept is originally due to 
Puig (see \cite{Pg-jalg}), and one version of his definition is given in 
Section 1 (Definition \ref{sat.Frob.}).  One motivating example is the 
fusion system of finite group $G$ with $S\in\sylp{G}$: the category 
$\calf_S(G)$ whose objects are the subgroups of $S$ and whose morphisms are 
those group homomorphisms which are conjugation by elements of $G$.

For $S\in\sylp{G}$ as above, there is a second, closely related category 
which can be defined, and which supplies the ``link'' between $\calf_S(G)$ 
and the classifying space $BG$ of $G$.  A 
subgroup $P\le{}S$ is called \emph{$p$-centric in $G$} if 
$Z(P)\in\sylp{C_G(P)}$; equivalently, if $C_G(P)=Z(P)\times{}C'_G(P)$ for 
some (unique) subgroup $C'_G(P)$ of order prime to $p$.  Let $\call_S^c(G)$ 
(the \emph{centric linking system of $G$}) be the category whose objects 
are the subgroups of $S$ which are $p$-centric in $G$, and where for each 
pair of objects $P,Q$:
	\[ \Mor_{\call_S^c(G)}(P,Q) = \bigl\{ g\in{}G \,\big|\, \9gP\le{}Q 
	\bigr\} \big/ C'_G(P)~. \]
Such categories were originally defined by Puig in \cite{Pg-mem}.  

To explain the significance of linking systems from a topologist's point of 
view, we must first define the \emph{geometric realization} of an arbitrary 
small category $\calc$.  This is a space $|\calc|$ built up of one vertex 
(point) for each object in $\calc$, one edge for each nonidentity morphism 
(with endpoints attached to the vertices corresponding to its source and 
target), one 2-simplex (triangle) for each commutative triangle in $\call$, 
etc.  (See, e.g., \cite[\S\,III.2.1--2]{AKO} for more details.)  By a 
theorem of Broto, Levi, and Oliver \cite[Proposition 1.1]{BLO1}, for any 
$G$ and $S$ as above, the space $|\call_S^c(G)|$, after $p$-completion in 
the sense of Bousfield and Kan, is homotopy equivalent to the $p$-completed 
classifying space $BG\pcom$ of $G$.  Furthermore, many of the homotopy 
theoretic properties of the space $BG\pcom$, such as its self homotopy 
equivalences, can be determined combinatorially by the properties (such as 
automorphisms) of the finite category $\call_S^c(G)$ \cite[Theorems B \& 
C]{BLO1}.  

Abstract \emph{centric linking systems} associated to a fusion system were 
defined in \cite{BLO2} (see Definition \ref{L-cat}).  
One of the motivations in \cite{BLO2} for defining these categories  
was that it provides a way to associate a classifying space to a 
saturated fusion system.  More precisely, if $\call$ is a centric linking 
system associated to a saturated fusion system $\calf$, then we regard 
the $p$-completion $|\call|\pcom$ of its geometric realization as a 
classifying space for $\calf$.  This is motivated by the equivalence 
$|\call_S^c(G)|\pcom\simeq{}BG\pcom$ noted above.  To give one example 
of the role played by these classifying spaces, if $\call'$ is another 
centric linking system, associated to a fusion system $\calf'$, and the 
classifying spaces $|\call|\pcom$ and $|\call'|\pcom$ are homotopy 
equivalent, then $\call\cong\call'$ and $\calf\cong\calf'$.  We refer to 
\cite[Theorem A]{BLO2} for more details and discussion.

It is unclear from the definition whether there is a centric linking system 
associated to any given saturated fusion system, and if so, whether it is 
unique.  Even when working with fusion systems of finite groups, which  
always have a canonical associated linking system, there is no simple 
reason why two groups with isomorphic fusion systems need have isomorphic 
linking systems, and hence equivalent $p$-completed classifying spaces. 
This question --- whether $\calf_S(G)\cong\calf_T(H)$ implies 
$\call_S^c(G)\cong\call_T^c(H)$ and hence $BG\pcom\simeq{}BH\pcom$ --- was 
originally posed by Martino and Priddy, and was what first got this author 
interested in the subject.

The main theorem of Chermak described in this paper is the following.

\begin{Th}[{Chermak \cite{Chermak}}] \label{ThA}
Each saturated fusion system has an associated centric linking system, 
which is unique up to isomorphism.
\end{Th}

\begin{proof} This follows immediately from Theorem \ref{lim*=0} in this 
paper, together with \cite[Proposition 3.1]{BLO2}. 
\end{proof}

In particular, this provides a new proof of the Martino-Priddy conjecture, 
which was originally proven in \cite{limz-odd,limz} using the 
classification of finite simple groups.  Chermak's theorem is much more 
general, but it also (indirectly) uses the classification in its proof.

Theorem \ref{ThA} is proven by Chermak by directly and systematically 
constructing the linking system, and by directly constructing an 
isomorphism between two given linking systems.  The proof given here 
follows the same basic outline, but uses as its main tool the obstruction 
theory which had been developed in \cite[Proposition 3.1]{BLO2} for dealing 
with this problem.  So if this approach is shorter, it is only because we 
are able to profit by the results of \cite[\S\,3]{BLO2}, and also by other 
techniques which have been developed more recently for computing these 
obstruction groups. 

By \cite[Proposition 4.6]{BLO3}, there is a bijective correspondence 
between centric linking systems associated to a given saturated fusion 
system $\calf$ up to isomorphism, and homotopy classes of rigidifications 
of the homotopy functor $\calo(\calf^c)\Right2{}\hotop$ which sends $P$ to 
$BP$.  (See Definition \ref{d:O(F)} for the definition of 
$\calo(\calf^c)$.)  Furthermore, if $\call$ corresponds to a 
rigidification $\til{B}$, then $|\call|$ is homotopy equivalent to the 
homotopy direct limit of $\til{B}$.  
Thus another consequence of Theorem \ref{ThA} is:

\begin{Th} \label{ThA2}
For each saturated fusion system $\calf$, there is a functor 
	\[ \til{B}\: \calo(\calf^c) \Right4{} \spaces, \]
together with a choice of homotopy equivalences $\til{B}(P)\simeq{}BP$ for 
each object $P$, such that for each 
$[\varphi]\in\Mor_{\calo(\calf^c)}(P,Q)$, the composite 
	\[ BP \simeq \til{B}(P) \Right4{\til{B}([\varphi])} \til{B}(Q) 
	\simeq BQ \]
is homotopic to $B\varphi$. Furthermore, $\til{B}$ is unique up to 
homotopy equivalence of functors, and $\hclim{}(\til{B})\pcom$ is the 
(unique) classifying space for $\calf$.  
\end{Th}

We also want to compare ``outer automorphism groups'' of fusion systems, 
linking systems, and their classifying spaces.  When $\calf$ is a saturated 
fusion system over a $p$-group $S$, set
	\[ \Aut(S,\calf)=\{\alpha\in\Aut(S)\,|\,\9\alpha\calf=\calf\} 
	\quad\textup{and}\quad 
	\Out(S,\calf)=\Aut(S,\calf)/\Aut_\calf(S)~. \]
Here, for $\alpha\in\Aut(S)$, $\9\alpha\calf$ is the fusion system over $S$ 
for which $\Hom_{\9\alpha\calf}(P,Q)=
\alpha\circ\homf(\alpha^{-1}(P),\alpha^{-1}(Q))\circ\alpha^{-1}$.  Thus 
$\Aut(S,\calf)$ is the group of ``fusion preserving'' automorphisms of $S$.

When $\call$ is a centric linking system associated to $\calf$, then for 
each object $P$ of $\call$, there is a ``distinguished monomorphism'' 
$\delta_P\:P\Right2{}\Aut_\call(P)$ (Definition \ref{L-cat}).  An 
automorphism $\alpha$ of $\call$ (a bijective functor from $\call$ to 
itself) is called \emph{isotypical} if it permutes the images of the 
distinguished monomorphisms; i.e., if 
$\alpha(\delta_P(P))=\delta_{\alpha(P)}(\alpha(P))$ for each $P$.  We 
denote by $\Out\typ(\call)$ the group of isotypical automorphisms of 
$\call$ modulo natural transformations of functors.  See also 
\cite[\S\,2.2]{AOV} or \cite[Lemma III.4.9]{AKO} for an alternative 
description of this group.  

By \cite[Theorem D]{BLO2}, $\Out\typ(\call)\cong\Out(|\call|\pcom)$, where 
$\Out(|\call|\pcom)$ is the group of homotopy classes of self homotopy 
equivalences of the space $|\call|\pcom$.  This is one reason for the 
importance of this particular group of (outer) automorphisms of $\call$.  
Another reason is the role played by $\Out\typ(\call)$ in the definition of 
a \emph{tame fusion system} in \cite[\S\,2.2]{AOV}.  

The other main consequence of the results in this paper is the following.

\begin{Th} \label{ThB}
For each saturated fusion system $\calf$ over a $p$-group $S$ with 
associated centric linking system $\call$, the natural homomorphism 
	\[ \Out\typ(\call) \Right4{\mu_\call} \Out(S,\calf) \]
induced by restriction to $\delta_S(S)\cong{}S$ is surjective, and is an 
isomorphism if $p$ is odd.
\end{Th}

\begin{proof} By \cite[III.5.12]{AKO}, 
$\Ker(\mu_\call)\cong\higherlim{}1(\calz_\calf)$, and $\mu_\call$ is onto 
whenever $\higherlim{}2(\calz_\calf)=0$.  (This was shown in \cite[Theorem 
E]{BLO1} when $\call$ is the linking system of a finite group.)  So the 
result follows from Theorem \ref{lim*=0} in this paper. 
\end{proof}

I would like to thank Assaf Libman for very carefully reading this 
manuscript, pointing out a couple gaps in the arguments, and making other 
suggestions for improving it.  And, of course, I very much want to thank 
Andy Chermak for solving this problem, which has taken up so much of my 
time for the past ten years.


\newsect{Notation and background}

We first briefly recall the definitions of saturated fusion systems and 
centric linking systems.  For any group $G$ and any pair of subgroups 
$H,K\le{}G$, set 
	\[ \Hom_G(H,K) = \bigl\{ c_g=(x\mapsto gxg^{-1}) \,\big|\, 
	g\in{}G,~ \9gH\le{}K \bigr\} \subseteq \Hom(H,K). \]
A \emph{fusion system} $\calf$ over a finite $p$-group $S$ is a
category whose objects are the subgroups of $S$, and whose 
morphism sets $\homf(P,Q)$ satisfy the following two conditions:
\begin{itemize}
\item $\Hom_S(P,Q)\subseteq\homf(P,Q)\subseteq\Inj(P,Q)$ for all $P,Q\le S$.

\item For each $\varphi\in\homf(P,Q)$, 
$\varphi^{-1}\in\homf(\varphi(P),P)$.
\end{itemize}
Two subgroups $P,P'\le{}S$ are called \emph{$\calf$-conjugate} if they are 
isomorphic in the category $\calf$.  Let $P^\calf$ denote the set of 
subgroups $\calf$-conjugate to $P$.

The following is the definition of a saturated fusion system first 
formulated in \cite{BLO2}.  Other (equivalent) definitions, including the 
original one by Puig, are discussed and compared in \cite[\S\S\,I.2 \& 
I.9]{AKO}.

\begin{Defi} \label{sat.Frob.}  
\renewcommand{\labelenumi}{\hskip-4pt$\bullet$}
Let $\calf$ be a fusion system over a $p$-group $S$.  
\begin{itemize} 
\item A subgroup $P\le{}S$ is \emph{fully centralized in $\calf$} if 
$|C_S(P)|\ge|C_S(Q)|$ for all $Q\in{}P^\calf$.
\item A subgroup $P\le{}S$ is \emph{fully normalized in $\calf$} if 
$|N_S(P)|\ge|N_S(Q)|$ for all $Q\in{}P^\calf$.
\item A subgroup $P\le{}S$ is \emph{$\calf$-centric} if 
$C_S(Q)\le{}Q$ for all $Q\in{}P^\calf$.
\item The fusion system $\calf$ is \emph{saturated} if the following 
two conditions hold:
\begin{enumerate}[\rm(I) ]
\item For all $P\le{}S$ which is fully normalized in $\calf$, $P$ is fully 
centralized in $\calf$ and $\Aut_S(P)\in\sylp{\autf(P)}$.
\item If $P\le{}S$ and $\varphi\in\homf(P,S)$ are such that $\varphi(P)$ is 
fully centralized, and if we set
	$$ N_\varphi = \{ g\in{}N_S(P) \,|\, \varphi c_g\varphi^{-1} \in 
	\Aut_S(\varphi(P)) \}, $$
then there is $\widebar{\varphi}\in\homf(N_\varphi,S)$ such that 
$\widebar{\varphi}|_P=\varphi$.
\end{enumerate}
\end{itemize}
\end{Defi}

For any fusion system $\calf$ over $S$, let $\calf^c\subseteq\calf$ be 
the full subcategory whose objects are the $\calf$-centric subgroups of 
$S$, and also let $\calf^c$ denote the set of $\calf$-centric subgroups of 
$S$.

\begin{Defi}[{\cite{BLO2}}]  \label{L-cat}
Let $\calf$ be a fusion system over the $p$-group $S$.  A \emph{centric 
linking system} associated to $\calf$ is a category $\call$ with 
$\Ob(\call)=\calf^c$, together with a functor $\pi\:\call\Right2{}\calf^c$
and \emph{distinguished monomorphisms} $P\Right1{\delta_P}\Aut_{\call}(P)$ 
for each $P\in\Ob(\call)$, which satisfy the following conditions.
\begin{enumerate}[\rm(A) ]
\item  $\pi$ is the identity on objects and is surjective on morphisms.
For each $P,Q\in\calf^c$, $\delta_P(Z(P))$ acts freely on 
$\Mor_{\call}(P,Q)$ by composition, and $\pi$ induces a bijection
	$$ \Mor_{\call}(P,Q)/\delta_P(Z(P)) \Right5{\cong} \homf(P,Q). $$

\item  For each $g\in{}P\in\calf^c$, $\pi$ sends 
$\delta_P(g)\in\Aut_{\call}(P)$ to $c_g\in\Aut_{\calf}(P)$.

\item  For each $P,Q\in\calf^c$, $\psi\in\Mor_{\call}(P,Q)$, and 
$g\in{}P$, $\psi\circ\delta_P(g)=\delta_Q(\pi(\psi)(g))\circ\psi$ in 
$\Mor_\call(P,Q)$.
\end{enumerate}
\end{Defi}

We next fix some notation for sets of subgroups of 
a given group.  For any group $G$, let $\cals(G)$ be the set of 
subgroups of $G$.  If $H\le{}G$ is any subgroup, set 
	\[ \II(H,G) = \{K\in\cals(G) \,|\, K\ge H\}. \]

\begin{Defi} 
Let $\calf$ be a saturated fusion system over a finite $p$-group $S$.  An 
\emph{interval} of subgroups of $S$ is a subset $\calr\subseteq\cals(S)$ 
such that $P<Q<R$ and $P,R\in\calr$ imply $Q\in\calr$.  An interval is 
\emph{$\calf$-invariant} if it is invariant under $\calf$-conjugacy.
\end{Defi}

Thus, for example, an $\calf$-invariant interval $\calr\subseteq\cals(S)$ 
is closed under overgroups if and only if $S\in\calr$.  Each 
$\calf$-invariant interval has the form $\calr{\sminus}\calr_0$ for some 
pair of $\calf$-invariant intervals $\calr_0\subseteq\calr$ which are 
closed under overgroups.

We next recall the obstruction theory to the existence and uniqueness of 
linking systems.

\begin{Defi} \label{d:O(F)}
Let $\calf$ be a saturated fusion system over a finite $p$-group $S$.  
\begin{enuma} 
\item Let $\calo(\calf^c)$ be the \emph{centric orbit category} of $\calf$:  
\ $\Ob(\calo(\calf^c))=\calf^c$, and 
	\[ \Mor_{\calo(\calf^c)}(P,Q) = \Inn(Q){\backslash}\homf(P,Q). \]

\item Let $\calz_\calf\:\calo(\calf^c)\op\Right2{}\Ab$ be the functor which 
sends $P$ to $Z(P)=C_S(P)$.  If $\varphi\in\homf(P,Q)$, and $[\varphi]$ 
denotes its class in $\Mor(\calo(\calf^c))$, then 
$\calz_\calf([\varphi])=\varphi^{-1}$ as a homomorphism from $Z(Q)=C_S(Q)$ 
to $Z(P)=C_S(P)$.

\item For any $\calf$-invariant interval $\calr\subseteq\calf^c$, let 
$\calz_\calf^\calr$ be the subquotient functor of $\calz_\calf$ where 
$\calz_\calf^\calr(P)=Z(P)$ if $P\in\calr$ and $\calz_\calf^\calr(P)=0$ 
otherwise.  

\item For each $\calf$-invariant interval $\calr\subseteq\calf^c$, we write 
for short
	\[ \LL\calf\calr 
	= \higherlim{\calo(\calf^c)}*(\calz_\calf^\calr)~; \]
i.e., the higher derived functors of the inverse limit of 
$\calz_\calf^\calr$.
\end{enuma}
\end{Defi}

We refer to \cite[\S\,III.5.1]{AKO} for more discussion of the functors 
$\higherlim{}*(-)$.

Thus $\calz_\calf=\calz_\calf^{\calf^c}$, and 
$\higherlim{}*(\calz_\calf)=\LL\calf{\calf^c}$.  By \cite[Proposition 
3.1]{BLO2}, the obstruction to the existence of a centric linking system 
associated to $\calf$ lies in $\LL[3]\calf{\calf^c}$, and the obstruction 
to uniqueness lies in $\LL[2]\calf{\calf^c}$.  

For any $\calf$ and any $\calf$-invariant interval $\calr$, 
$\calz_\calf^\calr$ is a quotient functor of $\calz_\calf$ if $S\in\calr$ 
(if $\calr$ is closed under overgroups).  If $\calr_0\subseteq\calr$ are 
both $\calf$-invariant intervals, and $P\in\calr_0$ and 
$Q\in\calr{\sminus}\calr_0$ implies $P\ngeq{}Q$, then 
$\calz_\calf^{\calr_0}$ is a subfunctor of $\calz_\calf^\calr$. 

\begin{Lem} \label{fixpt}
Fix a finite group $\Gamma$ with Sylow subgroup $S\in\sylp{\Gamma}$, and 
set $\calf=\calf_S(\Gamma)$.  Let $\calq\subseteq\calf^c$ be an 
$\calf$-invariant interval such that $S\in\calq$ (i.e., $\calq$ is closed 
under overgroups).  
\begin{enuma} 
\item 
Let $F\:\calo(\calf^c)\op\Right3{}\Ab$ be a functor such that $F(P)=0$ 
for each $P\in\calf^c{\sminus}\calq$.  
Let $\calo(\calf_\calq)\subseteq\calo(\calf^c)$ be the full subcategory 
with object set $\calq$. Then 
	\[ \higherlim{\calo(\calf^c)}*(F) \cong 
	\higherlim{\calo(\calf_\calq)}*(F|_{\calo(\calf_\calq)})~. \]

\item Assume $\calq=\II(Y,S)$ for some $p$-subgroup $Y\nsg\Gamma$ 
such that $C_\Gamma(Y)\le{}Y$.  Then 
	\beq \LL[k]\calf{\calq} \defeq 
	\higherlim{\calo(\calf^c)}k(\calz_\calf^\calq) 
	\cong \begin{cases} Z(\Gamma) & \textup{if $k=0$} \\ 
	0 & \textup{if $k>0$~.} \end{cases} \eeq
\end{enuma}
\end{Lem}

\begin{proof} \textbf{(a) }  Set $\calc=\calo(\calf^c)$ and 
$\calc_0=\calo(\calf_\calq)$ for short.  There is no morphism in $\calc$ 
from any object of $\calc_0$ to any object not in $\calc_0$.  Hence for any 
functor $F\:\calc\op\Right2{}\Ab$ such that $F(P)=0$ for each 
$P\notin\Ob(\calc_0)$, the two chain complexes $C^*(\calc;F)$ and 
$C^*(\calc_0;F|_{\calc_0})$ are isomorphic (see, e.g., 
\cite[\S\,III.5.1]{AKO}).  So 
$\higherlim{}*(F)\cong\higherlim{}*(F|_{\calc_0})$ in this situation, and 
this proves (a).  Alternatively, (a) follows upon showing that any 
$\calc_0$-injective resolution of $F|_{\calc_0}$ can be extended to an 
$\calc$-injective resolution of $F$ by assigning to all functors the value 
zero on objects not in $\calc_0$.

\smallskip

\noindent\textbf{(b) } To simplify notation, set $\4H=H/Y$ for each 
$H\in\II(Y,\Gamma)$, and $\4g=gY\in\4\Gamma$ for each 
$g\in\Gamma$.  
Let 
$\calo_{\4S}(\4\Gamma)$ be the ``orbit category'' of $\4\Gamma$:  the 
category whose objects are the subgroups of $\4S$, and where for 
$P,Q\in\calq$,
	\[ \Mor_{\calo_{\4S}(\4\Gamma)}(\4P,\4Q) = \4Q \big\backslash 
	\bigl\{ g\in{}\4\Gamma\,\big|\, \9g\4P\le \4Q \bigr\}~. \]
There is an isomorphism of categories 
$\Psi\:\calo(\calf_\calq)\Right2{\cong}\calo_{\4S}(\4\Gamma)$ which sends 
$P\in\calq$ to $\4P=P/Y$ and sends $[c_g]\in\Mor_{\calo(\calf_\calq)}(P,Q)$ 
to $\4Q\4g$.  Then $\calz_\calf^\calq\circ\Psi^{-1}$ sends $\4P$ to 
$Z(P)=C_{Z(Y)}(\4P)$.  Hence for $k\ge0$,
	\[ \higherlim{\calo(\calf^c)}k(\calz_\calf^\calq) \cong 
	\higherlim{\calo(\calf_\calq)}k
		(\calz_\calf^\calq|_{\calo(\calf_\calq)}) \cong
	\higherlim{\calo_{\4S}(\4\Gamma)}k
		(\calz_\calf^\calq\circ\Psi^{-1}) \cong 
	\begin{cases} C_{Z(Y)}(\4\Gamma)=Z(\Gamma) & \textup{if $k=0$}\\
	0 & \textup{if $k>0$} \end{cases} \]
where the first isomorphism holds by (a), and the last by a theorem of 
Jackowski and McClure \cite[Proposition 5.14]{JM}. We refer to 
\cite[Proposition 5.2]{JMO} for more details on the last isomorphism.
\end{proof}

More tools for working with these groups come from the long exact sequence 
of derived functors induced by a short exact sequence of functors.

\begin{Lem} \label{ex.seq.}
Let $\calf$ be a saturated fusion system over a finite $p$-group $S$.  Let 
$\calq$ and $\calr$ be $\calf$-invariant intervals such that
\begin{enumerate}[\rm(i) ] 
\item $\calq\cap\calr=\emptyset$,
\item $\calq\cup\calr$ is an interval, and
\item $Q\in\calq$, $R\in\calr$ implies $Q\nleq{}R$.
\end{enumerate}
Then $\calz_\calf^\calr$ is a subfunctor of $\calz_\calf^{\calq\cup\calr}$, 
$\calz_\calf^{\calq\cup\calr}\big/\calz_\calf^\calr\cong\calz_\calf^\calq$, 
and there is a long exact sequence
\begin{small} 
	\begin{multline*} 
	0 \Right2{} \LL[0]\calf\calr \Right2{} 
	\LL[0]\calf{\calq{\cup}\calr} \Right2{} \LL[0]\calf\calq \Right2{} 
	\cdots \\
	\Right2{} \LL[k-1]\calf\calq \Right2{} \LL[k]\calf\calr 
	\Right2{} \LL[k]\calf{\calq{\cup}\calr} \Right2{} \LL[k]\calf\calq 
	\Right2{} \cdots .
	\end{multline*}
\end{small}
In particular, the following hold.
\begin{enuma} 
\item If $\LL[k]\calf\calr\cong\LL[k]\calf\calq=0$ for some $k\ge0$, then 
$\LL[k]\calf{\calq\cup\calr}=0$.

\item Assume $\calf=\calf_S(\Gamma)$, where $S\in\sylp{\Gamma}$, and there is a 
normal $p$-subgroup $Y\nsg{}\Gamma$ such that $C_\Gamma(Y)\le{}Y$ and 
$\calq\cup\calr=\II(Y,S)$.  Then for each $k\ge2$, 
$\LL[k-1]\calf\calq\cong\LL[k]\calf\calr$.  Also, there is a short exact 
sequence
	\[ 1 \Right2{} C_{Z(Y)}(\Gamma) \Right4{} C_{Z(Y)}(\Gamma^*) 
	\Right4{} \LL[1]\calf\calr \Right2{} 1, \]
where $\Gamma^*=\Gen{g\in{}\Gamma\,\big|\,\9gP\in\calq \textup{ for some } 
P\in\calq}$.  
\end{enuma}
\end{Lem}

\begin{proof} Condition (iii) implies that $\calz_\calf^\calr$ is a 
subfunctor of $\calz_\calf^{\calq\cup\calr}$, and it is then immediate from 
the definitions (and (i) and (ii)) that 
$\calz_\calf^{\calq\cup\calr}\big/\calz_\calf^\calr\cong\calz_\calf^\calq$.  
The long exact sequence is induced by this short exact sequence of 
functors and the snake lemma.  Point (a) now follows immediately.  

Under the hypotheses in (b), by Lemma \ref{fixpt}(b), 
$\LL[k]\calf{\calq\cup\calr}=0$ for $k>0$ and 
$\LL[0]\calf{\calq\cup\calr}\cong Z(\Gamma)=C_{Z(Y)}(\Gamma)$.  The first 
statement in (b) thus follows immediately from the long exact sequence, and 
the second since $\LL[0]\calf\calq\cong{}C_{Z(Y)}(\Gamma^*)$ (by 
definition of inverse limits).  
\end{proof}

We next consider some tools for making computations in the groups 
$\higherlim{}*(-)$ for functors on orbit categories.  

\begin{Defi} \label{d:Lambda}
Fix a finite group $G$ and a $\Z[G]$-module $M$.  Let $\calo_p(G)$ be the 
category whose objects are the $p$-subgroups of $G$, and where 
$\Mor_{\calo_p(G)}(P,Q)=Q{\backslash}\{g\in{}G\,|\,\9gP\le{}Q\}$.  Define a 
functor $F_M\:\calo_p(G)\op\Right2{}\Ab$ by setting 
	\[ F_M(P)=\begin{cases} M & \textup{if $P=1$}\\
	0 & \textup{if $P\ne1$~.} 
	\end{cases} \]
Here, $F_M(1)=M$ has the given action of $\Aut_{\calo_p(G)}(1)=G$.  
Set
	$$ \Lambda^*(G;M) = \higherlim{\calo_p(G)}*(F_M). $$
\end{Defi}

These groups $\Lambda^*(G;M)$ provide a means of computing higher 
limits of functors on orbit categories which vanish except on one conjugacy 
class.

\begin{Prop}[{\cite[Proposition 3.2]{BLO2}}] \label{Lambda-red}
Let $\calf$ be a saturated fusion system over a $p$-group $S$. Let 
	\[ F\:\orb(\calf^c)\op\Right5{}\zploc\mod \] 
be any functor which vanishes except on the isomorphism class of 
some subgroup $Q\in\calf^c$.  Then
	\[ \higherlim{\calo(\calf^c)}*(F) \cong
	\Lambda^*(\outf(Q);F(Q)). \]
\end{Prop}

Upon combining Proposition \ref{Lambda-red} with the exact sequences of 
Lemma \ref{ex.seq.}, we get the following corollary.

\begin{Cor} \label{c:lim*=0}
Let $\calf$ be a saturated fusion system over a $p$-group $S$, and let 
$\calr\subseteq\calf^c$ be an $\calf$-invariant interval.  Assume, for some 
$k\ge0$, that $\Lambda^k(\outf(P);Z(P))=0$ for each $P\in\calr$.  Then 
$\LL[k]\calf\calr=0$.
\end{Cor}

What makes these groups $\Lambda^*(-;-)$ so useful is that they vanish in 
many cases, as described by the following proposition.

\begin{Prop}[{\cite[Proposition 6.1(i,ii,iii,iv)]{JMO}}] \label{Lambda}
The following hold for each finite group $G$ and each $\zploc[G]$-module 
$M$.
\begin{enuma}  
\item If $p{\nmid}|G|$, then 
$\Lambda^i(G;M)=\begin{cases}  M^G & \textup{if $i=0$}\\ 0 
& \textup{if $i>0$.} \end{cases}$

\item Let $H=C_G(M)$ be the kernel of the $G$-action on $M$.  Then 
$\Lambda^*(G;M)\cong\Lambda^*(G/H;M)$ if $p{\nmid}|H|$, and 
$\Lambda^*(G;M)=0$ if $p\big||H|$.

\item If $O_p(G)\ne1$, then $\Lambda^*(G;M)=0$.  

\item If $M_0\le{}M$ is a $\zploc[G]$-submodule, then there is an exact 
sequence
	\begin{multline*} 
	0 \Right2{} \Lambda^0(G;M_0) \Right2{} \Lambda^0(G;M) \Right2{} 
	\Lambda^0(G;M/M_0) \Right2{} \cdots \\
	\cdots \Right2{} \Lambda^{n-1}(G;M/M_0) \Right2{} \Lambda^n(G;M_0) 
	\Right2{} \Lambda^n(G;M) \Right2{} \cdots.
	\end{multline*}
\end{enuma}
\end{Prop}

The next lemma allows us in certain cases to replace the orbit category for 
one fusion system by that for a smaller one.  For any saturated fusion 
system $\calf$ over $S$ and any $Q\le{}S$, the \emph{normalizer fusion 
system} $N_\calf(Q)$ is defined as a fusion system over $N_S(Q)$ (cf. 
\cite[Definition I.5.3]{AKO}).  If $Q$ is fully normalized, then 
$N_\calf(Q)$ is always saturated (cf. \cite[Theorem I.5.5]{AKO}).  

\begin{Lem} \label{NF(Q)}
Let $\calf$ be a saturated fusion system over a $p$-group $S$, fix a 
subgroup $Q\in\calf^c$ which is fully normalized in $\calf$, 
and set $\cale=N_\calf(Q)$.  Set $\cale^\bullet=\calf^c\cap\cale^c$, a 
full subcategory of $\cale^c$, and let 
$\calo(\cale^\bullet)\subseteq\calo(\cale^c)$ be its orbit category. Define
	$$ \calt = \bigl\{ P\le S \,\big|\, 
	Q\nsg P, \textup{ and $R\in{}Q^\calf$, $R\nsg{}P$ implies $R=Q$}
	\bigr\}\,. $$
Let $F\:\calo(\calf^c)\op\Right2{}\zploc\mod$ be any functor which 
vanishes except on subgroups $\calf$-conjugate to subgroups in $\calt$, set 
$F_0=F|_{\calo(\cale^\bullet)}$, and let $F_1\:\calo(\cale^c)\op\Right2{}\Ab$ 
be such that $F_1|_{\calo(\cale^\bullet)}=F_0$ and $F_1(P)=0$ for all 
$P\in\cale^c{\sminus}\calf^c$.  Then 
restriction to $\cale^\bullet$ induces isomorphisms 
	\beqq \higherlim{\calo(\calf^c)}*(F) \RIGHT6{R}{\cong} 
	\higherlim{\calo(\cale^\bullet)}*(F_0) \LEFT6{R_1}{\cong}
	\higherlim{\calo(\cale^c)}*(F_1)~.
	\label{e:NFQ} \eeqq
\end{Lem}

\begin{proof} Since $R_1$ is an isomorphism by Lemma \ref{fixpt}(a), we 
only need to show that $R$ is an isomorphism.  If $F'\subseteq{}F$ is a 
pair of functors from $\calo(\calf^c)\op$ to $\zploc\mod$, and the lemma 
holds for $F'$ and for $F/F'$, then it also holds for $F$ by the 5-lemma 
(and since $R$ and $R_1$ are both natural in $F$ and preserve short exact 
sequences of functors).  It thus suffices to prove that the maps in 
\eqref{e:NFQ} are isomorphisms when $F$ vanishes except on the 
$\calf$-conjugacy class of one subgroup in $\calt$.  

Fix $P\in\calt$, and assume $F(R)=0$ for all $R\notin{}P^\calf$.  Thus 
$Q\nsg{}P$ by definition of $\calt$.  If $\varphi\in\homf(P,S)$ is such 
that $Q\nsg\varphi(P)$, then $\varphi^{-1}(Q)\nsg{}P$, so $\varphi(Q)=Q$ 
($\varphi\in\Aut_\cale(Q)$) since $P\in\calt$.  Thus 
$\Out_\cale(P)=\outf(P)$, and $Q\nsg{}P^*\in{}P^\calf$ implies 
$P^*\in{}P^\cale$. This yields the following diagram:
	\beq \vcenter{\xymatrix@C=30pt@R=40pt{
	\higherlim{\calo(\calf^c)}*(F) \ar@<1.2ex>[r]^-{R} 
	\ar[dr]^-{\Phi^*}_-{\cong} & 
	\higherlim{\calo(\cale^\bullet)}*(F_0) &
	\higherlim{\calo(\cale^c)}*(F_1) 
	\ar@<-1.2ex>[l]_-{R_1}^-{\cong} \ar[dl]_-{\Phi_1^*}^-{\cong} \\
	& \Lambda^*(\outf(P);F(P)) 
	\rlap{\,,}
	}} \eeq
where $\Phi^*$ and $\Phi_1^*$ are the 
isomorphisms of Proposition \ref{Lambda-red}. The commutativity of the 
diagram follows from the precise description of $\Phi^*$ and $\Phi^*_1$ in 
\cite[Proposition 3.2]{BLO2}.  Thus $R$ is an isomorphism.
\end{proof}

The following lemma can also be stated and proven as a result about 
extending automorphisms from a linking system to a group \cite[Lemma 
4.17]{Chermak}.

\begin{Lem} [{\cite[4.17]{Chermak}}] \label{ch4.17}
Fix a pair of finite groups $H\nsg{}G$, together with $S\in\sylp{G}$ and 
$T=S\cap{}H\in\sylp{H}$.  Set $\calf=\calf_S(G)$ and $\cale=\calf_T(H)$.  
Assume $Y\le{}T$ is such that $Y\nsg{}G$ and $C_G(Y)\le{}Y$.  Let $\calq$ 
be an $\calf$-invariant interval in $\II(Y,S)$ such that 
$S\in\calq$, and such that $Q\in\calq$ implies $H\cap{}Q\in\calq$.  Set 
$\calq_0=\{Q\in\calq\,|\,Q\le{}H\}$.  Then restriction induces an injective 
homomorphism 
	\[ \LL[1]\calf\calq \Right4{R} \LL[1]\cale{\calq_0}. \]
\end{Lem}

\begin{proof} Since $\cale^c$ need not be contained in $\calf^c$, we must 
first check that there is a well defined ``restriction'' homomorphism.  
Set $\cale^\bullet=\cale^c\cap\calf^c$:  a full subcategory 
of $\cale^c$.  Since the functor $\calz_\cale^{\calq_0}$ vanishes on all 
subgroups in $\cale^c$ not in $\calq_0\subseteq\cale^\bullet$, the higher 
limits are the same whether taken over $\calo(\cale^\bullet)$ or 
$\calo(\cale^c)$ (Lemma \ref{fixpt}(a)).  Thus $R$ is defined as the 
restriction map to $\higherlim{}1(\calz_\cale^{\calq_0}|_{\cale^\bullet}) 
\cong\LL[1]\cale{\calq_0}$. 

We work with the bar resolutions for $\calo(\calf^c)$ and 
$\calo(\cale^\bullet)$, using the notation of \cite[\S\,III.5.1]{AKO}.  
Fix a cocycle $\eta\in{}Z^1(\calo(\calf^c);\calz_\calf^\calq)$ such that 
$[\eta]\in\Ker(R)$.  Thus $\eta$ is a function from 
$\Mor(\calo(\calf^c))$ to $Z(Y)$ which sends the class $[\varphi]$ of 
$\varphi\in\Hom_G(P,Q)$ to an element of $Z(P)$ if $P\in\calq$, and to $1$ 
if $P\notin\calq$. We can assume, after adding an appropriate coboundary, 
that $\eta\bigl(\Mor(\calo(\cale^\bullet))\bigr)=1$.  

Define $\widehat{\eta}\in{}Z^1(N_G(T)/T;Z(T))$ to be the restriction of 
$\eta$ to $\Aut_{\calo(\calf^c)}(T)=N_G(T)/T$.  For $g\in{}N_G(T)$, let 
$\bar{g}$ be its class in $N_G(T)/T$. Set 
$\gamma=\eta([\incl_T^S])\in{}Z(T)$, so $d\gamma\in{}Z^1(N_G(T)/T;Z(T))$ is 
the cocycle $d\gamma(\bar{g})=\gamma^g\cdot\gamma^{-1}$. For each 
$g\in{}S$, $[\incl_T^S]\circ[c_g]=[\incl_T^S]$ in $\calo(\calf^c)$, so 
$\gamma^g\cdot\eta([c_g])=\gamma$, and thus 
$\widehat{\eta}(\bar{g})=\eta([c_g])=(d\gamma(\bar{g}))^{-1}$.  In 
other words, $\widehat{\eta}|_{S/T}$ is a coboundary, and since 
$S/T\in\sylp{N_G(T)/T}$, $[\widehat{\eta}]=1\in{}H^1(N_G(T)/T;Z(T))$ (cf. 
\cite[Theorem XII.10.1]{CE}).  Hence there is $\beta\in{}Z(T)$ such that 
$\widehat{\eta}=d\beta$.  Since $\eta([c_h])=1$ for all $h\in{}N_H(T)$, 
$[\beta,h]=1$ for all $h\in{}N_H(T)$, and thus $\beta\in{}Z(N_H(T))$.  

Let $G^*\le{}G$ be the subgroup generated by all $g\in{}G$ such that for 
some $Q\in\calq$, $\9gQ\in\calq$.  Define $H^*\le{}H$ similarly.  Since 
$S\le{}N_G(T)\le{}G^*$ and $N_H(T)\le{}H^*$, $S\in\sylp{G^*}$, 
$T\in\sylp{H^*}$, and $HG^*\ge{}HN_G(T)=G$ by the Frattini argument (Lemma 
\ref{NP>P}(b)).  If $g=ha$ where $h\in{}H$, $a\in{}N_G(T)$, and 
$\9gQ\in\calq$ for some $Q\in\calq$, then $\9a(Q\cap{}H)$ 
and $\9g(Q\cap{}H)=\9gQ\cap{}H$ are both in $\calq_0$, and thus 
$h\in{}H^*$.  Since $N_G(T)$ normalizes $H^*$, this shows that 
$G^*=H^*N_G(T)$.  So $G^*\cap{}H=(H^*N_G(T))\cap{}H=H^*N_H(T)=H^*$.  In 
particular, $H^*\nsg{}G^*$ and $G^*/H^*\cong{}G/H$. 

For each $\varphi\in\Hom_H(P,Q)$ (where $Y\le{}P,Q\le{}T$), and each 
$g\in{}N_G(T)$, set $\9g\varphi=c_g\varphi{}c_g^{-1}\in\Hom_H(\9gP,\9gQ)$.  
Since $\eta([\varphi])=\eta([\9g\varphi])=1$, we have 
$\varphi^{-1}(\widehat{\eta}(\bar{g}))=\widehat{\eta}(\bar{g})$.  Thus for 
each $g\in{}N_G(T)$, $\widehat{\eta}(\bar{g})=\beta^g\beta^{-1}$ is 
invariant under the action of $H^*$; i.e., $\beta^g\beta^{-1}\in{}Z(H^*)$.  
So the class $[\beta]\in{}Z(N_H(T))/Z(H^*)$ is fixed under the action of 
$N_G(T)$ on this quotient.


Since $p{\nmid}[H^*:N_H(T)]$, and since $N_G(T)$ normalizes $H^*$ and 
$N_H(T)$, the inclusion of $Z(H^*)=C_{Z(Y)}(H^*)$ into 
$Z(N_H(T))=C_{Z(Y)}(N_H(T))$ is $N_G(T)$-equivariantly split by the trace 
homomorphism for the actions of $H^*\ge{}N_H(T)$ on $Z(Y)$.  So the fixed 
subgroup for the $N_G(T)$-action on the quotient group $Z(N_H(T))/Z(H^*)$ 
is $Z(N_G(T))/Z(G^*)$.  Thus $\beta\in{}Z(N_G(T))Z(H^*)$, and we can assume 
$\beta\in{}Z(H^*)$ without changing $d\beta=\widehat{\eta}$.  

Define a $0$-cochain 
$\widehat{\beta}\in{}C^0(\calo(\calf^c);\calz_\calf^\calq)$ by setting 
$\widehat{\beta}(P)=\beta$ if $P\in\calq_0$ and 
$\widehat{\beta}(P)=1$ otherwise.  Then 
$\eta([\varphi])=d\widehat{\beta}([\varphi])$ for all 
$\varphi\in\Mor(\cale^\bullet)$ (since both vanish) and also for all 
$\varphi\in\Aut_G(T)$.  Since $G=HN_G(T)$, each morphism in 
$\calf$ between subgroups of $T$ is the composite of a morphism in $\cale$ 
and the restriction of a morphism in $\autf(T)$.  Hence 
$\eta([\varphi])=d\widehat{\beta}([\varphi])$ for all such morphisms 
$\varphi$ (since $\eta$ and $d\widehat{\beta}$ are both cocycles).  Upon 
replacing $\eta$ by $\eta(d\widehat{\beta})^{-1}$, we can assume $\eta$ 
vanishes on all morphisms in $\calf$ between subgroups of $T$.

For each $P\in\calq$, set $P_0=P\cap{}T$ and let $i_P\in\Hom_G(P_0,P)$ be 
the inclusion.  Then $\eta([i_P])\in{}Z(P_0)$ (and $\eta([i_P])=1$ if 
$P\notin\calq$).  For each $g\in{}P$, the relation $[i_P]=[i_P]\circ[c_g]$ 
in $\calo(\calf^c)$ (where $[c_g]\in\Aut_{\calo(\calf^c)}(P_0)$) implies 
that $\eta([i_P])$ is $c_g$-invariant.  Thus $\eta([i_P])\in{}Z(P)$.  Let 
$\rho\in{}C^0(\calo(\calf^c);\calz_\calf^\calq)$ be the 0-cochain 
$\rho(P)=\eta([i_P])$ when $P\in\calq$ and $\rho(P)=1$ if 
$P\in\calf^c{\sminus}\calq$.  Thus $\rho(P)=1$ if $P\le{}T$ by the initial 
assumptions on $\eta$. Then $d\rho([i_P])=\eta([i_P])$ for each $P$, and 
$d\rho(\varphi)=1=\eta(\varphi)$ for each $\varphi$ between subgroups of 
$T$.  For each $\varphi\in\Hom_G(P,Q)$ in $\calf^c$, let 
$\varphi_0\in\Hom_G(P_0,Q_0)$ be its restriction; the relation 
$[\varphi]\circ[i_Q]=[i_P]\circ[\varphi_0]$ in $\calo(\calf^c)$ implies 
that $\eta([\varphi])=d\rho([\varphi])$ since this holds for $[\varphi_0]$ 
and the inclusions.  Thus $\eta=d\rho$, and so $[\eta]=1$ in 
$\LL[1]\calf\calq$.  \end{proof}

We end the section by recalling a few elementary results about finite 
groups.

\begin{Lem} \label{NP>P}
\begin{enuma} 
\item If $Q>P$ are $p$-groups for some prime $p$, then $N_Q(P)>P$.
\item \textup{(Frattini argument) } If $H\nsg{}G$ are finite groups and 
$T\in\sylp{H}$, then $G=HN_G(T)$.
\end{enuma}
\end{Lem}

\begin{proof} See, for example, \cite[Theorems 2.1.6 \& 2.2.7]{Sz1}.
\end{proof}

\begin{Lem} \label{CG(V)}
Let $G$ be a finite group such that $O_p(G)=1$, and assume $G$ acts 
faithfully on an abelian $p$-group $D$.  Then $G$ acts faithfully on 
$\Omega_1(D)$.
\end{Lem}

\begin{proof} The subgroup $C_G(\Omega_1(D))$ is a normal $p$-subgroup of 
$G$ (cf. \cite[Theorem 5.2.4]{Gorenstein}), and hence is contained in 
$O_p(G)=1$. 
\end{proof}

\newsect{The Thompson subgroup and offenders}

The proof of the main theorem is centered around the Thompson subgroup of 
a $p$-group, and the $FF$-offenders for an action of a group on an abelian 
$p$-group.  We first fix the terminology and notation which will be used.

\begin{Defi} \label{d:offenders}
\begin{enuma} 
\item For any $p$-group $S$, set $d(S)=\sup\bigl\{|A|\,\big|\,A\le{}S \textup{ 
abelian}\bigr\}$, let $\cala(S)$ be the set of all abelian 
subgroups of $S$ of order $d(S)$, and set $J(S)=\gen{\cala(S)}$.

\item Let $G$ be a finite group which acts faithfully on the abelian 
$p$-group $D$.  A \emph{best offender} in $G$ on $D$ is an abelian subgroup 
$A\le{}G$ such that $|A||C_D(A)|\ge|B||C_D(B)|$ for each $B\le{}A$.  (In 
particular, $|A||C_D(A)|\ge|D|$.)  Let $\cala_D(G)$ be the set of best 
offenders in $G$ on $D$, and set $J_D(G)=\gen{\cala_D(G)}$.

\item Let $\Gamma$ be a finite group, and let $D\nsg{}\Gamma$ be a normal 
abelian $p$-subgroup.  Let $J(\Gamma,D)\le{}\Gamma$ be the subgroup such 
that $J(\Gamma,D)/C_\Gamma(D)=J_D(\Gamma/C_\Gamma(D))$.  
\end{enuma}
\end{Defi}

Note, in the situation of point (c) above, that
	\[ D \le C_\Gamma(D)\le J(\Gamma,D)\le \Gamma 
	\qquad\textup{and}\qquad
	J(J(\Gamma,D),D)=J(\Gamma,D). \]

The relation between the Thompson subgroup $J(-)$ and best offenders is 
described by the next lemma and corollary.

\begin{Lem} \label{l:offenders}
\begin{enuma} 
\item Assume $G$ acts faithfully on a finite abelian $p$-group $D$. If $A$ 
is a best offender in $G$ on $D$, and $U$ is an $A$-invariant subgroup of $D$, 
then $A/C_A(U)$ is a best offender in $N_G(U)/C_G(U)$ on $U$.

\item Let $S$ be a finite $p$-group, let $D\nsg{}S$ be a normal abelian 
subgroup, and set $G=S/C_S(D)$.  Assume $A\in\cala(S)$.  Then the image of 
$A$ in $G$ is a best offender on $D$. 

\end{enuma}
\end{Lem}

\begin{proof} We give here the standard proofs.

\smallskip

\textbf{(a) } Set $\4A=A/C_A(U)$ for short.  For each
$\4B=B/C_A(U)\le\4A$,
	\[ |C_U(B)||C_D(A)| = |C_U(B)C_D(A)||C_U(B)\cap{}C_D(A)|
	\le |C_D(B)||C_U(A)|. \]
Also, $|B||C_D(B)|\le|A||C_D(A)|$ since $A$ is a best offender on $D$, and hence 
	\[ |\4B||C_U(\4B)| = \frac{|B||C_U(B)|}{|C_A(U)|}
	\le \frac{|B||C_D(B)|}{|C_D(A)|}\cdot\frac{|C_U(A)|}{|C_A(U)|}
	\le |A|\cdot\frac{|C_U(A)|}{|C_A(U)|}
	= |\4A||C_U(\4A)|. \]
Thus $\4A$ is a best offender on $U$.

\smallskip

\textbf{(b) } Set $\4A=A/C_A(D)$, identified with the image of $A$ in $G$.  
Fix some $\4B=B/C_A(D)\le\4A$, and set $B^*=C_D(B)B$.  This is an abelian 
group since $D$ and $B$ are abelian, and hence $|B^*|\le|A|$ since 
$A\in\cala(S)$.  Since $B\cap{}C_D(B)\le{}C_D(A)$, 
	\[ |\4B||C_D(\4B)| 
	= \frac{|B||C_D(B)|}{|C_A(D)|} 
	= \frac{|B^*||B\cap{}C_D(B)|}{|C_A(D)|} 
	\le \frac{|A||C_D(A)|}{|C_A(D)|} = |\4A||C_D(\4A)|. \]
Since this holds for all $\4B\le\4A$, $\4A$ is a best offender on $D$.
\end{proof}

The following corollary reinterprets Lemma \ref{l:offenders} in terms of 
the groups $J(\Gamma,D)$ defined above.

\begin{Cor} \label{c:offenders}
Let $\Gamma$ be a finite group, and let $D\nsg{}\Gamma$ be a normal abelian 
$p$-subgroup.
\begin{enuma} 
\item If $U\le{}D$ is also normal in $\Gamma$, then 
$J(\Gamma,U)\ge J(\Gamma,D)$.
\item If $\Gamma$ is a $p$-group, then $J(\Gamma)\le J(\Gamma,D)$.
\end{enuma}
\end{Cor}

An action of a group $G$ on a group $D$ is \emph{quadratic} if 
$[G,[G,D]]=1$.  If $D$ is abelian and $G$ acts faithfully, then a 
\emph{quadratic best offender} in $G$ on $D$ is an abelian subgroup 
$A\le{}G$ which is a best offender and whose action is quadratic.

\begin{Lem} \label{quad.act} 
Let $G$ be a finite group which acts faithfully on an elementary abelian 
$p$-group $V$.  If the action of $G$ on $V$ is quadratic, then $G$ is also 
an elementary abelian $p$-group. 
\end{Lem}

\begin{proof} We write $V$ additively for convenience; thus $[g,v]=gv-v$ 
for $g\in{}G$ and $v\in{}V$.  By an easy calculation, and since the action 
is quadratic, $[gh,v]=[g,v]+[h,v]$ for each $g,h\in{}G$ and $v\in{}V$.  
Thus $g\mapsto(v\mapsto[g,v])$ is a homomorphism from $G$ to 
the additive group $\End(V)$, and is injective since the action is 
faithful.  Since $\End(V)$ is an elementary abelian $p$-group, so is $G$.
\end{proof}

We will also need the following form of Timmesfeld's replacement theorem.

\begin{Thm} \label{Timmesfeld}
Let $A$ and $V$ be abelian $p$-groups.  Assume $A$ acts on $V$, and 
is a best offender on $V$.  Then there is $1\ne{}B\le{}A$ such 
that $B$ is a quadratic best offender on $V$.  More precisely, we can take 
$B=C_A([A,V])\ne1$, in which case $|A||C_V(A)|=|B||C_V(B)|$ and 
$C_V(B)=[A,V]+C_V(A)<V$.
\end{Thm}

\begin{proof} We follow the proof given by Chermak in 
\cite[\S\,1]{Chermak-quad}.  Set $m=|A||C_V(A)|$.  Since $A$ is a best 
offender, 
	\beqq m = \sup\bigl\{ |B||C_U(B)| \,\big|\, B\le{}A,~ U\le{}V 
	\bigr\}.  \label{e:2.5a} \eeqq
For each $U\le{}V$, consider the set 
	\[ \calm_U = \bigl\{ B \le A \,\big|\, |B||C_U(B)|=m \bigr\}~. \]
By the maximality of $m$ in \eqref{e:2.5a}, 
	\beqq B\in\calm_U \quad\implies\quad |C_U(B)|=|C_V(B)| 
	\quad\implies\quad C_V(B)\le U~. \label{e:2.5b} \eeqq

\smallskip

\noindent\textbf{Step 1: } (Thompson's replacement theorem)  
For each $x\in{}V$, set 
	\[ V_x=[A,x] \defeq \Gen{[a,x]=ax-x \,\big|\, a\in{}A}
	\qquad\textup{and}\qquad A_x=C_A(V_x). \]
Note that $V_x$ is $A$-invariant.  We will show that 
	\beqq |A_x||C_V(A_x)|=|A||C_V(A)|=m \qquad\textup{and}\qquad
	C_V(A_x) = V_x + C_V(A). \label{e:2.5d} \eeqq

Define $\Phi\:A\Right2{}V_x$ (a map of sets) by setting 
$\Phi(a)=[a,x]=ax-x$ for each $a\in{}A$.  We first claim that $\Phi$ 
induces an injective map of sets
	\[ \phi\: A/A_x \Right4{} V_x/C_{V_x}(A) \]
between these quotient groups.  Since $A$ is abelian, 
$[a,[b,x]]=abx-bx-ax+x=[b,[a,x]]$ for all $a,b\in{}A$. Hence for all 
$g,h\in{}A$, 
	\begin{align*} 
	\Phi(g)-\Phi(h)= & gx-hx \in C_{V_x}(A) 
	~\iff~ h([h^{-1}g,x])\in C_{V_x}(A) \\
	&~\iff~ 1 = [A,[h^{-1}g,x]] = [h^{-1}g,[A,x]] = [h^{-1}g,V_x] \\
	&~\iff~ h^{-1}g \in C_A(V_x)=A_x~. 
	\end{align*}
Thus $\phi$ is well defined and injective.

Now, 
	\beqq |V_x||C_V(A)| = |C_{V_x}(A)||V_x+C_V(A)| \le 
	|C_{V_x}(A)||C_V(A_x)|, \label{e:2.5c} \eeqq
since $V_x\le{}C_V(A_x)$ by definition of $A_x$.  Together with the 
injectivity of $\phi$, this implies that
	\[ \frac{|A|}{|A_x|} \le \frac{|V_x|}{|C_{V_x}(A)|} \le 
	\frac{|C_V(A_x)|}{|C_V(A)|}, \]
and so $m=|A||C_V(A)|\le|A_x||C_V(A_x)|$.  The opposite inequality holds by 
\eqref{e:2.5a}, so $A_x\in\calm_V$ and the inequality in \eqref{e:2.5c} is 
an equality.  Thus $|V_x+C_V(A)|=|C_V(A_x)|$, finishing the proof of 
\eqref{e:2.5d}.  

\smallskip

\noindent\textbf{Step 2: } Assume, for some $U\le{}V$, that 
$B_0,B_1\in\calm_U$.  Then $m=|B_0||C_U(B_0)|\ge|B_0B_1||C_U(B_0B_1)|$ by 
\eqref{e:2.5a}, and hence 
	\[ \frac{|B_1|}{|B_0\cap{}B_1|} = \frac{|B_0B_1|}{|B_0|} 
	\le \frac{|C_U(B_0)|}{|C_U(B_0B_1)|} = 
	\frac{|C_U(B_0)+C_U(B_1)|}{|C_U(B_1)|} \le
	\frac{|C_U(B_0\cap{}B_1)|}{|C_U(B_1)|}. \]
So $m=|B_1||C_U(B_1)|\le|B_0\cap{}B_1||C_U(B_0\cap{}B_1)|$ with 
equality by \eqref{e:2.5a} again, and we conclude that 
$B_0\cap{}B_1\in\calm_U$.
	
\smallskip

\noindent\textbf{Step 3: } Set $B=C_A([A,V])$ and $U=[A,V]+C_V(A)$.  For 
each $x\in{}V$, \eqref{e:2.5d} implies that $C_V(A_x)=[A,x]+C_V(A)\le{}U$ 
and $A_x\in\calm_U$.  Hence $B=\bigcap_{x\in{}V}A_x\in\calm_U$ by Step 2, 
so $B$ is a best offender on $V$, and is quadratic since $[B,[A,V]]=1$ by 
definition.  Also, $C_V(B)\le{}U$ by \eqref{e:2.5b}.  Since 
$U=[A,V]+C_V(A)\le{}C_V(B)$ by definition, we conclude that $U=C_V(B)$.

If $U=V$, then $V=[A,V]\oplus{}W$ is an $A$-invariant splitting for some 
$W\le{}C_V(A)$.  But this would imply $[A,V]=[A,[A,V]]+[A,W]=[A,[A,V]]$, 
which is impossible since $[A,X]<X$ for any finite $p$-group $X$ on which 
$A$ acts.  We conclude that $U=C_V(B)<V$, and hence that $B\ne1$. 
\end{proof}


\newsect{Proof of the main theorem}

The following terminology will be very useful when carrying out the reduction 
procedures used in this section.

\begin{Defi}[{\cite[6.3]{Chermak}}] 
A \emph{general setup} is a triple $(\Gamma,S,Y)$, where $\Gamma$ is a 
finite group, $S\in\sylp{\Gamma}$, $Y\nsg{}\Gamma$ is a normal 
$p$-subgroup, and $C_\Gamma(Y)\le{}Y$ ($Y$ is centric in $\Gamma$).  A 
\emph{reduced setup} is a general setup $(\Gamma,S,Y)$ such that 
$Y=O_p(\Gamma)$, $C_S(Z(Y))=Y$, and $O_p(\Gamma/C_\Gamma(Z(Y)))=1$. 
\end{Defi}

The next proposition, which will be shown in Section 4, is the 
key technical result needed to prove the main theorem.  Its proof uses the 
classification by Meierfrankenfeld and Stellmacher \cite{MS2} of 
$FF$-offenders, and through that depends on the classification of finite 
simple groups.

\begin{Prop}[{Compare \cite[6.11]{Chermak}}] \label{ch8.10}
Let $(\Gamma,S,Y)$ be a reduced setup, set $D=Z(Y)$, and 
assume $\Gamma/C_\Gamma(D)$ is generated by quadratic best offenders on $D$.  Set 
$\calf=\calf_S(\Gamma)$, and let $\calr\subseteq\calf^c$ be the set of all 
$R\ge{}Y$ such that $J(R,D)=Y$.  Then $\LL[2]\calf\calr=0$ if $p=2$, and 
$\LL[1]\calf\calr=0$ if $p$ is odd.
\end{Prop}

Since this distinction between the cases where $p=2$ or where $p$ is odd 
occurs throughout this section and the next, it will be convenient to 
define
	\[ k(p) = \begin{cases} 2 & \textup{if $p=2$} \\
	1 & \textup{if $p$ is an odd prime.}
	\end{cases} \]
Thus under the hypotheses of Proposition \ref{ch8.10}, we claim that 
$\LL[k(p)]\calf\calr=0$.

Proposition \ref{ch8.10} seems very restricted in scope, but it can be 
generalized to the following situation.

\begin{Prop}[{Compare \cite[6.12]{Chermak}}] \label{ch8.11}
Let $(\Gamma,S,Y)$ be a general setup, and assume Proposition \ref{ch8.10} holds 
for all reduced setups $(\Gamma^*,S^*,Y^*)$ with $|\Gamma^*|\le|\Gamma|$.  Set 
$\calf=\calf_S(\Gamma)$ and $D=Z(Y)$.  Let $\calr\subseteq\II(Y,S)$ 
be an $\calf$-invariant interval such that 
for each $Q\in\II(Y,S)$, $Q\in\calr$ if and only if 
$J(Q,D)\in\calr$.  
Then $\LL[k]\calf\calr=0$ for all $k\ge k(p)$.
\end{Prop}

\begin{proof} Assume the proposition is false.  Let $(\Gamma,S,Y,\calr,k)$ be a 
counterexample for which the 4-tuple $(k,|\Gamma|,|\Gamma/Y|,|\calr|)$ is the 
smallest possible under the lexicographical ordering.  


We will show in Step 1 that $\calr=\{P\le{}S\,|\,J(P,D)=Y\}$, in Step 2 
that $k=k(p)$, in Step 3 that $(\Gamma,S,Y)$ is a reduced setup, and in 
Step 4 that $\Gamma/C_\Gamma(D)$ is generated by quadratic best offenders 
on $D$.  The result then follows from Proposition \ref{ch8.10}.

\smallskip

\noindent\textbf{Step 1: } Let $R_0\in\calr$ be a minimal element of 
$\calr$ which is fully normalized in $\calf$.  Since $J(R_0,D)\in\calr$ by 
assumption (and $J(R_0,D)\le{}R_0$), $J(R_0,D)=R_0$.  Let $\calr_0$ be the 
set of all $R\in\calr$ such that $J(R,D)$ is $\calf$-conjugate to $R_0$, 
and set $\calq_0=\calr{\sminus}\calr_0$.  Then $\calr_0$ and $\calq_0$ are 
both $\calf$-invariant intervals, and satisfy the conditions $Q\in\calr_0$ 
($Q\in\calq_0$) if and only if $J(Q,D)\in\calr_0$ ($J(Q,D)\in\calq_0$).  
Since $\LL[k]\calf\calr\ne0$, Lemma \ref{ex.seq.}(a) implies 
$\LL[k]\calf{\calr_0}\ne0$ or $\LL[k]\calf{\calq_0}\ne0$. Hence 
$\calq_0=\emptyset$ and $\calr=\calr_0$ by the minimality assumption on 
$|\calr|$ (and since $\calr_0\ne\emptyset$). 

Set $\Gamma_1=N_\Gamma(R_0)$, $S_1=N_S(R_0)$, 
$\calf_1=N_\calf(R_0)=\calf_{S_1}(\Gamma_1)$ (see \cite[Proposition 
I.5.4]{AKO}), and $\calr_1=\{R\in\calr\,|\,J(R,D)=R_0\}$.  Every subgroup 
in $\calr$ is $\calf$-conjugate to a subgroup in $\calr_1$.  Also, for 
$P\in\calr_1$, if $R\in{}R_0{}^\calf$ and $R\nsg{}P$, then 
$J(P,D)\ge{}J(R,D)=R$ implies $R=R_0$.  The hypotheses of Lemma \ref{NF(Q)} 
are thus satisfied, and so 
$\LL[k]{\calf_1}{\calr_1}\cong\LL[k]\calf\calr\ne0$.  Thus 
$(\Gamma_1,S_1,Y,\calr_1,k)$ is another a counterexample to the 
proposition.  By the minimality assumption, $\Gamma_1=\Gamma$, and thus 
$R_0\nsg{}\Gamma$.

We have now shown that there is a $p$-subgroup $R_0\nsg{}\Gamma$ 
such that $\calr=\{R\le{}S\,|\,J(R,D)=R_0\}$.  Set $Y_1=R_0\ge{}Y$ and 
$D_1=Z(Y_1)\le{}D$.  By Corollary \ref{c:offenders}(a), for each $R\le{}S$ 
such that $R\ge{}R_0$ and $R\notin\calr$, $J(R,D_1)\ge{}J(R,D)\notin\calr$, 
and hence $J(R,D_1)\notin\calr$.  Thus $(\Gamma,S,Y_1,\calr,k)$ is a 
counterexample to the proposition, and so $Y=Y_1=R_0$ by the minimality 
assumption on $|\Gamma/Y|$.  We conclude that $\calr=\{R\le{}S\,|\,J(R,D)=Y\}$. 

\smallskip

\noindent\textbf{Step 2: } Let $\calq$ be the set of all overgroups of $Y$ 
in $S$ which are not in $\calr$.  Equivalently, 
$\calq=\{Q\le{}S\,|\,J(Q,D)>Y\}$.  If $k\ge2$, then 
$\LL[k-1]\calf\calq\cong\LL[k]\calf\calr\ne0$ by Lemma \ref{ex.seq.}(b).  
Since $k$ was assumed to be the smallest degree $\ge{}k(p)$ for which the 
proposition is not true, we conclude that $k=k(p)$.

\smallskip

\noindent\textbf{Step 3: } Assume $(\Gamma,S,Y)$ is not a reduced setup.  
Let $K\nsg{}\Gamma$ be such that $K\ge{}C_\Gamma(D)$ and 
$K/C_\Gamma(D)=O_p(\Gamma/C_\Gamma(D))$, and set $Y_2=S\cap{}K\nsg{}S$.  
Then $Y_2>Y$, since either $Y_2\ge{}O_p(\Gamma)>Y$, or $Y_2\ge{}C_S(D)>Y$, 
or $p\big||K/C_\Gamma(D)|$ and hence $Y_2>C_S(D)\ge{}Y$.  Set 
$\Gamma_2=N_\Gamma(Y_2)$, and set $\calr_2=\{P\in\calr\,|\,P\ge{}Y_2\}$.  
Note that $S\in\sylp{\Gamma_2}$, and also that $\calr_2$ is an 
$\calf$-invariant interval since $Y_2$ is strongly closed in $S$ with 
respect to $\Gamma$.  Set $\calf_2=\calf_S(\Gamma_2)=N_\calf(Y_2)$ 
\cite[Proposition I.5.4]{AKO}.

Assume $P\in\calr{\sminus}\calr_2$.  Then $P\ngeq{}Y_2$, so $PY_2>P$, and 
hence $N_{PY_2}(P)>P$ (Lemma \ref{NP>P}(a)).  Set $G=\Out_\Gamma(P)$ and 
$G_0=\Out_K(P)$.  Then $G_0\nsg{}G$ since $K\nsg{}\Gamma$, and 
$C_{G_0}(Z(P))=\Out_{C_K(Z(P))}(P)\ge \Out_{C_\Gamma(D)}(P)$ since 
$K\ge{}C_\Gamma(D)$ and $Z(P)\le{}Z(Y)=D$.  Hence $G_0/C_{G_0}(Z(P))$ is a 
$p$-group since $K/C_\Gamma(D)$ is a $p$-group.  For any 
$g\in{}N_{PY_2}(P){\sminus}P$, $\Id\ne[c_g]\in\Out_K(P)=G_0$ since 
$Y_2\le{}K$ (and since $C_\Gamma(P)\le{}C_\Gamma(Y)\le{}Y\le{}P$).  Thus 
$\Out_K(P)=G_0\nsg{}G$ contains a nontrivial element of $p$-power order, 
and its action on $Z(P)$ factors through the $p$-group $G_0/C_{G_0}(Z(P))$.  
Proposition \ref{Lambda}(b,c) now implies that 
$\Lambda^*(\Out_\Gamma(P);Z(P))=0$.  

Since this holds for all $P\in\calr{\sminus}\calr_2$, 
$\LL\calf{\calr{\sminus}\calr_2}=0$ by Corollary \ref{c:lim*=0}.  Hence 
$\LL\calf{\calr_2}\cong\LL\calf\calr$ by the exact sequence in Lemma 
\ref{ex.seq.}.  Also, the hypotheses of Lemma \ref{NF(Q)} hold for the 
functor $\calz_\calf^{\calr_2}$ on $\calo(\calf^c)$ (with $Q=Y_2$) since 
$Y_2$ is strongly closed.  So 
$\LL\calf{\calr_2}\cong\LL{\calf_2}{\calr_2}$.  Since 
$\LL[k]\calf\calr\ne0$ by assumption, $\LL[k]{\calf_2}{\calr_2}\ne0$.  

Set $D_2=Z(Y_2)\le{}D$. For each 
$P\in\II(Y_2,S)$, 
	\beqq  P \ge J(P,D_2) \ge J(P,D) \ge C_P(D) \ge Y \label{e:3.2a} \eeqq
by Corollary \ref{c:offenders}(a) and by definition of $J(P,-)$. 
We must show that 
$P\in\calr_2$ if and only if $J(P,D_2)\in\calr_2$. If 
$P\in\calr_2\subseteq\calr$, then $J(P,D)\in\calr$ by assumption, so 
$J(P,D_2)\in\calr$ by \eqref{e:3.2a} since $\calr$ is an interval, and 
$J(P,D_2)\in\calr_2$ since $J(P,D_2)\ge{}C_P(D_2)\ge{}Y_2$.  If 
$P\notin\calr_2$, then $P\notin\calr$, so $J(P,D)\notin\calr$, and 
$J(P,D_2)\notin\calr$ (hence $J(P,D_2)\notin\calr_2$) by \eqref{e:3.2a} 
again and since $\calr$ is an interval containing $Y$.  

Thus $(\Gamma_2,S,Y_2,\calr_2,k)$ is a counterexample 
to the proposition.  So $\Gamma_2=\Gamma$ and $Y_2=Y$ by the minimality 
assumption, which contradicts the above claim that $Y_2>Y$.  We conclude 
that $(\Gamma,S,Y)$ is a reduced setup.

\smallskip

\noindent\textbf{Step 4: } It remains to prove that $\Gamma/C_\Gamma(D)$ is 
generated by quadratic best offenders on $D$; the result then 
follows from Proposition \ref{ch8.10}.  

Let $\Gamma_3\nsg{}\Gamma$ be such that 
$\Gamma_3\ge{}C_\Gamma(D)$ and $\Gamma_3/C_\Gamma(D)$ is generated by all 
quadratic best offenders on $D$.  If $\Gamma_3=\Gamma$ we are done, so 
assume $\Gamma_3<\Gamma$.  Set $S_3=\Gamma_3\cap{}S$ and 
$\calf_3=\calf_{S_3}(\Gamma_3)$.  Set 
	\[ \calq = \II(Y,S){\sminus}\calr,\quad
	\calq_3 = \calq\cap\II(Y,S_3), \quad\textup{and}\quad
	\calr_3 = \calr\cap\II(Y,S_3). \]
Since $\LL[k]\calf\calr\ne0$, $\calr\subsetneqq\II(Y,S)$ by Lemma 
\ref{fixpt}(b), and $\calq\ne\emptyset$.  
The proposition holds for $(\Gamma_3,S_3,Y,\calr_3,k)$ by the minimality 
assumption, and thus $\LL[k]{\calf_3}{\calr_3}=0$. 

For $Q\in\calq$, $J(Q,D)>Y$, so $Q/Y$ has nontrivial best offenders on $D$, 
hence has nontrivial quadratic best offenders on $D$ by Theorem 
\ref{Timmesfeld}, and thus $J(Q\cap\Gamma_3,D)>Y$.  So $Q\in\calq$ implies 
$Q\cap\Gamma_3\in\calq_3$ by Step 1. In particular, $S_3\in\calq_3$.  

If $k=2$ (i.e., if $p=2$), then $\LL[1]\calf\calq\cong\LL[2]\calf\calr\ne0$ 
and $\LL[1]{\calf_3}{\calq_3}\cong\LL[2]{\calf_3}{\calr_3}=0$ by Lemma 
\ref{ex.seq.}(b), which is impossible by Lemma \ref{ch4.17}.  

If $k=1$ (if $p$ is odd), set 
	\begin{align*} 
	\Gamma^*&=\Gen{g\in{}\Gamma\,\big|\,\9gP\in\calq 
	\textup{ for some } P\in\calq} \le \Gamma \\
	\Gamma_3^*&=\Gen{g\in{}\Gamma_3\,\big|\,\9gP\in\calq_3 
	\textup{ for some } P\in\calq_3} \le \Gamma_3~. 
	\end{align*}
Then $\Gamma_3^*\le\Gamma^*$ since $\Gamma_3\le\Gamma$ and 
$\calq_3\subseteq\calq$.  By Lemma \ref{ex.seq.}(b), there 
are exact sequences 
	\beqq \begin{split} 
	&1 \Right2{} C_{Z(Y)}(\Gamma^*) \Right4{} C_{Z(Y)}(\Gamma^*)  
	\Right4{} \LL[1]\calf\calq \ne1 \\
	&1 \Right2{} C_{Z(Y)}(\Gamma_3) \Right4{} C_{Z(Y)}(\Gamma_3^*) 
	\Right4{} \LL[1]{\calf_3}{\calq_3} =1.
	\end{split} \label{e:ww} \eeqq
Also, $\Gamma^*\Gamma_3\ge{}N_\Gamma(S_3)\Gamma_3=\Gamma$ since 
$S_3\in\calq_3$, where the equality follows from the Frattini argument 
(Lemma \ref{NP>P}(b)), so 
	\[ C_{Z(Y)}(\Gamma) = C_{Z(Y)}(\Gamma^*\Gamma_3) 
	= C_{Z(Y)}(\Gamma^*) \cap C_{Z(Y)}(\Gamma_3) . \]
But this is impossible, since $C_{Z(Y)}(\Gamma)<C_{Z(Y)}(\Gamma^*) \le 
C_{Z(Y)}(\Gamma_3^*)=C_{Z(Y)}(\Gamma_3)$ by the exactness in \eqref{e:ww}. 
\end{proof}

We now have the tools needed to prove the main vanishing result.

\begin{Thm} \label{lim*=0}
For each saturated fusion system $\calf$ over a $p$-group $S$, 
$\higherlim{\calq(\calf^c)}k(\calz_\calf)=0$ for all 
$k\ge2$, and for all $k\ge1$ if $p$ is odd.
\end{Thm}

\begin{proof} As in \cite[\S\,6]{Chermak}, we choose inductively 
subgroups $X_0,X_1,\ldots,X_N\in\calf^c$ and 
$\calf$-invariant intervals $\emptyset=\calq_{-1}\subseteq\calq_0\subseteq 
\cdots\subseteq\calq_N=\calf^c$ as follows.  
Assume $\calq_{n-1}$ has been defined ($n\ge0$), and 
$\calq_{n-1}\subsetneqq\calf^c$. Consider the sets of subgroups
	\begin{align*} 
	\calu_1 &= \bigl\{P\in\calf^c{\sminus}\calq_{n-1} \,\big|\, d(P) 
	\textup{ maximal} \bigr\} \\
	\calu_2 &= \bigl\{P\in\calu_1 \,\big|\, |J(P)| \textup{ maximal} 
	\bigr\} \\
	\calu_3 &= \bigl\{P\in\calu_2 \,\big|\, J(P)\in\calf^c \bigr\} \\
	\calu_4 &= \begin{cases} 
	\bigl\{P\in\calu_3 \,\big|\, |P| \textup{ minimal} \bigr\} & 
	\textup{if $\calu_3\ne\emptyset$} \\
	\bigl\{P\in\calu_2 \,\big|\, |P| \textup{ maximal} \bigr\} & 
	\textup{otherwise.}
	\end{cases}
	\end{align*}
(See Definition \ref{d:offenders}(a) for the definition of $d(P)$.)  Let 
$X_n$ be any subgroup in $\calu_4$ which is fully normalized in $\calf$.  
Since $\calu_4$ is invariant under $\calf$-conjugacy, this is always 
possible.

Let $\calq_n$ be the union of $\calq_{n-1}$ with the set of all overgroups 
of subgroups $\calf$-conjugate to $X_n$.  Set 
$\calr_n=\calq_n{\sminus}\calq_{n-1}$ for each 
$0\le{}n\le{}N$. By definition of $\calu_4$, $X_n=J(X_n)$ if 
$J(X_n)\in\calf^c$, while $\calr_n=X_n{}^\calf$ if $J(X_n)\notin\calf^c$.  
Note also that $X_0=J(S)$ and $\calr_0=\calq_0=\II(J(S),S)$. 

We will show, for each $n$, that 
	\beqq
	\textup{$\LL[k]\calf{\calr_n}=0$ for all $k\ge k(p)$.}
	\label{e:xx} \eeqq
Then by Lemma \ref{ex.seq.}(a), for all $k\ge{}k(p)$, 
$\LL[k]\calf{\calq_{n-1}}=0$ implies $\LL[k]\calf{\calq_n}=0$.  The theorem 
now follows by induction on $n$.

\smallskip

\noindent\textbf{Case 1: } Assume $n$ is such that $J(X_n)\notin\calf^c$, 
and hence that $\calr_n=X_n{}^\calf$.  Since 
$J(X_n)$ is centric in $X_n$ but not in $S$, $X_nC_S(J(X_n))>X_n$.  
Hence $N_{X_nC_S(J(X_n))}(X_n)>X_n$ (Lemma \ref{NP>P}(a)), so 
there is $g\in{}N_S(X_n){\sminus}X_n$ such that $[g,J(X_n)]=1$.  Then $g$ 
acts trivially on $Z(X_n)\le{}J(X_n)$, so the kernel of the 
$\outf(X_n)$-action on $Z(X_n)$ has order a multiple of $p$, and 
$\Lambda^*(\outf(X_n);Z(X_n))=0$ by Proposition \ref{Lambda}(b).  Hence 
\eqref{e:xx} holds by Proposition \ref{Lambda-red}.

\smallskip

\noindent\textbf{Case 2: } Assume $n$ is such that $J(X_n)\in\calf^c$, and 
hence $X_n=J(X_n)$ by definition of $\calu_4$.  By definition of $\calu_1$ 
and $\calu_2$, for each $P\ge{}X_n$ in $\calr_n$, $d(P)=d(X_n)$ and 
$J(P)=X_n$.  Hence 
	\beqq P\in\calr_n ~\implies~
	\textup{$J(P)$ is the unique subgroup of $P$ 
	$\calf$-conjugate to $X_n$.} \label{e:zz} \eeqq

Set $T=N_S(X_n)$ and $\cale=N_\calf(X_n)$.  Then $\cale$ is a saturated 
fusion system over $T$ (cf. \cite[Theorem I.5.5]{AKO}), and contains $X_n$ 
as normal centric subgroup.  Hence there is a \emph{model} for $\cale$ (cf. 
\cite[Theorem III.5.10]{AKO}):  a finite group $\Gamma$ such that 
$T\in\sylp{\Gamma}$, $X_n\nsg\Gamma$, $C_\Gamma(X_n)\le{}X_n$, and 
$\calf_T(\Gamma)=\cale$.  

Let $\calr$ be the set of all $P\in\calr_n$ such that $P\ge{}X_n$; thus 
$\calr=\calr_n\cap\cale$ by \eqref{e:zz}. Then 
$(\Gamma,T,X_n)$ is a general setup, and $\calr$ is an $\cale$-invariant 
interval containing $X_n$.  If $P\in\calr$ and $Y\le{}P$ is 
$\calf$-conjugate to $X_n$, then $Y=X_n$ by \eqref{e:zz}.  
The hypotheses of Lemma \ref{NF(Q)} thus hold, and hence 
	\beqq \LL\calf{\calr_n}\cong\LL\cale\calr. \label{e:3.4x} \eeqq

Set $D=Z(X_n)$.  We claim that for each $P\in\II(X_n,T)$,
	\beqq P\in\calr \quad \iff \quad J(P,D)\in\calr~. \label{e:yy} 
	\eeqq
Fix such a $P$. By Corollary \ref{c:offenders}(b), $J(P,D)\ge{}J(P)$, and 
$X_n\ge{}J(X_n,D)\ge{}J(X_n)=X_n$.  If $P\in\calr$, then $J(P,D)\in\calr$ 
since $X_n=J(X_n,D)\le J(P,D)\le P$ and $\calr$ is an interval.  If 
$P\notin\calr$, then $P\in\calq_{n-1}$, so $n\ge1$, and by definition of 
$\calu_1$ and $\calu_2$, either $d(P)>d(X_n)$, or $d(P)=d(X_n)$ and 
$J(P)>X_n$.  If $d(P)>d(X_n)$, then $d(J(P,D))=d(P)>d(X_n)$ since 
$J(P)\le{}J(P,D)\le{}P$, and $J(P,D)\notin\calr$ since $d(R)=d(X_n)$ for 
all $R\in\calr$.  If $J(P)>X_n$, then $J(P)\notin\calr$ since $J(R)=X_n$ 
for all $R\in\calr$, and hence $J(P,D)\notin\calr$ since $J(P,D)\ge{}J(P)$ 
and $\calr$ is an interval.  This proves \eqref{e:yy}.  

Thus by Proposition \ref{ch8.11}, $\LL[k]\cale\calr=0$ for all $k\ge k(p)$.  
Together with \eqref{e:3.4x}, this finishes the proof of \eqref{e:xx}, and 
hence of the theorem. 
\end{proof}


\newcommand{\ubm}{\textit{\textbf{m}}} 


\newcommand{\bfb}{{\textbf{\textit{b}}}}

\newsect{Proof of Proposition \ref{ch8.10}}

It remains to prove Proposition \ref{ch8.10}, which we restate here as:

\begin{Prop} \label{ch9.1}
Let $(\Gamma,S,Y)$ be a reduced setup, set $D=Z(Y)$, and 
assume $\Gamma/C_\Gamma(D)$ is generated by quadratic best offenders on $D$.  Set 
$\calf=\calf_S(\Gamma)$, and let $\calr\subseteq\calf^c$ be the set of all 
$R\ge{}Y$ such that $J(R,D)=Y$.  Then $\LL[k(p)]\calf\calr=0$.
\end{Prop}

It is in this section that we use the classification of offenders by 
Meierfrankenfeld and Stellmacher \cite{MS2}, and through that the 
classification of finite simple groups.  Any proof of Proposition 
\ref{rad.p-ch1} without using these results would imply a classification-free 
proof of Proposition \ref{ch9.1}, and hence of Theorem \ref{lim*=0}. 

In the following definition, when $H_1\le{}H_2\le\cdots\le{}H_k$ are 
subgroups of a group $G$, then $N_G(H_1,\ldots,H_k)$ denotes the 
intersection of their normalizers.

\begin{Defi} \label{d:rad.p-ch}
Let $G$ be a finite group.
\begin{enuma} 
\item A \emph{radical $p$-subgroup} of $G$ is a $p$-subgroup 
$P\le{}G$ such that $O_p(N_G(P))=P$; i.e., $O_p(N_G(P)/P)=1$.

\item A \emph{radical $p$-chain of length $k$} in $G$ is a sequence of 
$p$-subgroups $P_0<P_1<\cdots<P_k\le{}G$ such that $P_0$ is radical in $G$, 
$P_i$ is radical in $N_G(P_0,\ldots,P_{i-1})$ for each 
$i\ge1$, and $P_k\in\sylp{N_G(P_0,\ldots,P_{k-1})}$.  
\end{enuma}
\end{Defi}

The reason for defining this here is the following vanishing result, which 
involves only radical $p$-chains with $P_0=1$. 

\begin{Prop}[{\cite[Lemma III.5.27]{AKO} and \cite[Proposition 3.5]{limz}}] 
\label{Lambda-rad.ch}
Fix a finite group $G$, a finite $\F_p[G]$-module $M$, and $k\ge1$ such 
that $\Lambda^k(G;M)\ne0$.  Then there is a radical $p$-chain 
$1=P_0<P_1<\cdots<P_k$ of length $k$ such that $M$ contains a copy of the 
free module $\F_p[P_k]$.  
\end{Prop}

Other results similar to Proposition \ref{Lambda-rad.ch} were shown by 
Grodal, and follow from results in \cite[\S\,5]{Grodal}.

Note that since the trivial subgroup is a radical $p$-subgroup of $G$ only 
if $O_p(G)=1$, Proposition \ref{Lambda-rad.ch} includes the statement that 
$\Lambda^k(G;M)=0$ if $O_p(G)\ne1$ (Proposition \ref{Lambda}(c)). The 
reason for defining radical $p$-chains more generally here --- to also 
allow chains where $P_0\ne1$ --- will be seen in the statement of the 
following proposition and in the proof of Proposition \ref{ch9.1}.  

\begin{Prop} \label{rad.p-ch1}
Let $G$ be a nontrivial finite group with $O_p(G)=1$, and let $V$ be a 
faithful $\F_p[G]$-module.  Fix a set $\calu$ of nontrivial quadratic best 
offenders in $G$ on $V$ which is invariant under $G$-conjugacy, and assume 
$G=\gen{\calu}$.  Set $G_0=O^p(G)$ and $W=C_V(G_0)[G_0,V]/C_V(G_0)$.  
Set $k=k(p)$, and let 
$P_0<P_1<\cdots<P_k$ be a radical $p$-chain in $G$ such that $P_0$ does 
not contain any subgroup in $\calu$.  Then $C_W(P_0)$, with its induced 
action of $P_k/P_0$, does not contain a copy of the free module 
$\F_p[P_k/P_0]$.
\end{Prop}

\begin{proof} Quadratic offenders with faithful action are elementary 
abelian by Lemma \ref{quad.act}.  So the results in \cite{MS2} can be 
applied. In particular, $G=M=J_M(V)=J$ in the notation of \cite{MS2}, 
$O^p(G)=\gen{\calj}$ by \cite[Theorem 1(c)]{MS2}, and hence $W$ as defined 
here is the same as $W$ defined in \cite{MS2}.  

By \cite[Theorem 1(d,e)]{MS2}, $W$ is a semisimple $\F_p[G]$-module, and 
each $U\in\calu$ is a quadratic best offender on $W$.  Set 
$W=W_1\oplus\cdots\oplus{}W_m$, where each $W_i$ is $\F_p[G]$-irreducible.  
Set $K_i=C_G(W_i)$.  

Assume $C_W(P_0)$ contains a copy of $\F_p[P_k/P_0]$.  Then by Lemma 
\ref{l:A1}(a), there is $1\le j\le m$ such that $C_{W_j}(P_0)$ also 
contains a copy of $\F_p[P_k/P_0]$.  

Choose $S\in\sylp{G}$ such that $P_k\le{}S$, and set
$T=S\cap{}K_j\in\sylp{K_j}$.  Set $T_0=T\cap{}P_0$ and 
$T^*=N_T(P_k)$.  We claim that $T_0=T^*=T$.  Since 
$W_j$ contains $\F_p[P_k/P_0]$, $C_{P_0}(W_j)=C_{P_k}(W_j)$, and hence 
$T_0=P_0\cap{}T=P_k\cap{}T$.  If $T_0<T$, then $P_k\ngeq{}T$,
so $N_{P_kT}(P_k)>P_k$ by Lemma \ref{NP>P}(a), and 
$T^*=N_T(P_k)>T_0$.  For each $0\le i\le k$, 
	\[ [P_i,T^*] \le [P_k,T^*]\le{}P_k\cap{}T=T_0\le{}P_0 \le P_i, \]
so $P_i\nsg T^*$.  Thus $P_k<P_kT^*\le{}N_G(P_0,\ldots,P_{k-1})$, 
contradicting the assumption that $P_k\in\sylp{N_G(P_0,\ldots,P_{k-1})}$.  
We conclude that $T=T_0\le{}P_k$.  

Since $G=\gen{\calu}$ is generated by quadratic best offenders on $W$, 
$G/K_j$ is generated by quadratic best offenders on $W_j$  by Lemma 
\ref{l:offenders}(a).  So by Lemma \ref{rad.p-ch}, 
$P_kK_j/K_j\in\sylp{G/K_j}$.  Thus $P_k=S\in\sylp{G}$, since 
$P_k\ge{}T\in\sylp{K_j}$.  

Since $P_k\in\sylp{G}$ and $\calu$ is $G$-conjugacy invariant, there is 
$U\in\calu$ such that $U\le{}P_k$.  By assumption, $U\nleq{}P_0$. Set 
$U_0=U\cap{}P_0$, and set $p^r=|U/U_0|>1$.  Since 
$C_W(U_0)\ge{}C_W(P_0)$ contains a copy of the free 
module $\F_p[U/U_0]$ of rank $p^r$, and since 
$|U||C_W(U)|\ge|U_0||C_W(U_0)|$, 
	\[ p^r-1 \le \rk\bigl(C_W(U_0)/C_W(U)\bigr) \le 
	\log_p|U/U_0| = r. \]
Hence $p^r=|U/U_0|=2$; i.e., $p=2$ and $r=1$.  But $k=k(p)=2$ since $p=2$, so 
$|P_k/P_0|\ge4$, and $C_W(U_0)$ contains at least two copies of the 
free module $\F_2[U/U_0]$ since it contains a copy of $\F_2[P_k/P_0]$.  This 
implies $2\le\rk\bigl(C_W(U)/C_W(U_0)\bigr)\le{}r=1$.  Since this is impossible, 
we conclude that $C_W(P_0)$ does not contain a copy of $\F_p[P_k/P_0]$.  
\end{proof}

The following lemma was needed to prove Proposition \ref{rad.p-ch1}.

\begin{Lem} \label{rad.p-ch}
Let $G$ be a nontrivial finite group, let $W$ be a faithful, irreducible 
$\F_p[G]$-module, and assume $G$ is generated by its quadratic best 
offenders on $W$.  Let $P_0<P_1<\cdots<P_k$ be a radical $p$-chain in $G$ 
with $k\ge k(p)$.  Then either
\begin{enuma} 
\item $C_W(P_0)$, with its induced action of $P_k/P_0$, does not contain a 
copy of the free module $\F_p[P_k/P_0]$; or

\item $P_k\in\sylp{G}$.
\end{enuma}
\end{Lem}

\begin{proof} Quadratic offenders with faithful action are elementary 
abelian by Lemma \ref{quad.act}.  Also, $O_p(G)=1$, since 
$C_W(O_p(G))$ is a nontrivial $\F_p[G]$-submodule of $W$ and $G$ acts 
faithfully.  So the results in \cite{MS2} can be applied. Since the set 
$\cald$ defined in \cite[Theorem 2]{MS2} contains all quadratic offenders, 
$G=\gen{\cald}$.

By \cite[Theorem 2]{MS2}, either $G$ is a ``genuine'' group of Lie 
type in characteristic $p$; or $p=2$ and $G$ is one of the groups $3A_6$, 
$A_7$, $\Sigma_n$ or $A_n$ with the natural $\F_2[G]$-module, or 
$SO_{2m}^\pm(2^a)$ with the natural module.  We consider these cases 
individually.

\smallskip

\noindent\textbf{Case 1: } Assume $G$ is a genuine group of Lie type in 
characteristic $p$; in particular, $G$ is quasisimple.  The nontrivial 
radical $p$-subgroups of $G$ are well known:  by a theorem of Borel and 
Tits (see \cite[corollary]{BW} or \cite[13-5]{GL}), they are all conjugate 
to maximal normal unipotent subgroups in parabolic subgroups.  Hence the 
successive normalizers all contain Sylow $p$-subgroups of $G$, and the 
quotients are again groups of Lie type.  Since 
$P_k\in\sylp{N_G(P_0,\dots,P_{k-1})}$, $P_k\in\sylp{G}$ in this case.

\smallskip

\noindent\textbf{Case 2: } Assume $p=2$, and $G\cong3A_6$ or $A_7$.  Then 
the Sylow 2-subgroups of $G$ have order 8, the nontrivial radical 
$2$-subgroups have order 4 or 8, hence are normal in Sylow 2-subgroups, and 
thus $P_k\in\syl2{G}$ ($k\ge2$).

\smallskip

\noindent\textbf{Case 3: } Assume $p=2$, $G\cong\Sigma_m$ or $A_m$, and $W$ 
is a natural module for $G$. Thus $\rk(W)=m-1$ if $m$ is odd, while 
$\rk(W)=m-2$ if $m$ is even.  

Set $\ubm=\{1,2,\ldots,m\}$, with the action of $G$.  Fix a radical 
$p$-chain $P_0<P_1<\cdots<P_k$ with $k\ge2$, and assume $C_W(P_0)$ does 
contain a free submodule $\F_2[P_k/P_0]$.  Thus 
$\rk(C_W(P_0))\ge|P_k/P_0|\ge4$.  

If $P_0$ acts on $\ubm$ with more than one orbit, then by Lemma \ref{l:A4}, 
$|\ubm/P_0|\ge3$, and $C_W(P_0)$ is contained in $\F_2(\ubm/P_0)/\Delta$ 
where $\Delta\cong\F_2$ is generated by the sum of the basis elements.  Then 
$\F_2(\ubm/P_0)$ also contains a copy of the free module 
$\F_2[P_k/P_0]$, and hence $\ubm/P_0$ contains a free $(P_k/P_0)$-orbit by 
Lemma \ref{l:A1}(b).  But this is impossible by Lemma \ref{l:A3}.

Now assume $P_0$ acts transitively on $\ubm$.  We have 
$P_0=\widebar{P}_0\cap{}G$, where $\widebar{P}_0=O_2(N_{\Sigma_m}(P_0))$ is 
radical in $\Sigma_m$. By \cite[p. 7]{AF}, $\widebar{P}_0=E_{s_1}\wr 
E_{s_2}\wr\cdots\wr E_{s_n}\le\Sigma_m$ ($n\ge1$), where $E_s\cong C_2^s$ 
with a free action on a set of order $2^s$, and where $P_0$ acts on $\ubm$ 
via the canonical action of a wreath product ($m=2^{s_1+\ldots+s_n}$).  Let 
$Q\le{}P_0$ be the stabilizer of a point under the action on $\ubm$.  Then 
$P_0/\gen{\Fr(P_0),Q}\cong E_{s_n}$, so $C_W(P_0)\cong\F_2^{s_n}$ by Lemma 
\ref{l:A4}.  Also, $m\ge8$ since $4\big||P_k/P_0|\big||N_G(P_0)/P_0|$, so 
either $\widebar{P}_0\nleq{}A_m$ (if $s_1=1$), or 
$N_{\Sigma_m}(P_0)\le{}A_m$ (if $s_1\ge3$, or $s_1=2$ and $n\ge2$).  In 
both cases, $N_G(P_0)/P_0\cong 
N_{\Sigma_m}(\widebar{P}_0)/\widebar{P}_0\cong\prod_{i=1}^nGL_{s_i}(2)$ 
(\cite[p. 6]{AF}), and acts on $C_W(P_0)$ via projection to the factor 
$GL_{s_n}(\F_2)$.  As seen in Case 1, any radical $2$-chain in 
$GL_{s_n}(\F_2)$ must end with a Sylow 2-subgroup; i.e., $P_k/P_0\cong 
U\in\syl2{GL_{s_n}(\F_2)}$.  But the condition 
$|P_k/P_0|=2^{s_n(s_n-1)/2}\le\rk(C_W(P_0))=s_n$ then implies $s_n\le2$, 
which is impossible since $|P_k/P_0|\ge4$.  

\smallskip

\noindent\textbf{Case 4: } Now assume $p=2$ and $G\cong{}SO_{2m}^\pm(q)$, 
where $q=2^a$ ($a\ge1$), and $W$ is the natural $\F_2[G]$-module of rank 
$2am$.  Then $G\not\cong{}SO_4^+(2)\cong{}GL_2(2)\wr{}C_2$, since 
$G$ is not generated by best offenders in this case ($W$ 
splits in a unique way $W=W_1\oplus{}W_2$ where $W_1\perp{}W_2$ 
and $\rk(W_i)=2$, and all best offenders send each $W_i$ to itself).  
Also, $G\not\cong{}SO_4^-(2)\cong\Sigma_5$ (with the reduced permutation 
action on $W$) since this was already eliminated in Case 3.  Thus either 
$m\ge3$, or $m=2$ and $q\ge4$ ($a\ge2$).  Set $G_0=\Omega_{2m}^\pm(q)$, so 
$[G:G_0]=2$.  

For any nontrivial radical $2$-subgroup $1\ne{}P\le{}G$, $P\cap{}G_0$ is a 
radical $p$-subgroup of $G_0$, and hence is either trivial, or is a maximal 
normal unipotent subgroup in a parabolic subgroup.  If $P\cap{}G_0=1$ (and 
$P\ne1$), then $P=\gen{t}$ for some involution 
$t\in{}SO_{2m}^\pm(q){\sminus}\Omega_{2m}^\pm(q)$.  By \cite[8.10]{AS}, $t$ 
has type $b_\ell$ (in the notation of \cite[\S\,7--8]{AS}), where 
$\ell=\rk([t,W])$ is odd.  By \cite[8.7]{AS}, $O_2(N_G(t))$ contains all 
transvections $u\in{}G$ with $[u,W]\le[t,W]$, and hence $P=\gen{t}$ is 
radical only if $t$ is itself a transvection ($\ell=1$).  In that case, by 
\cite[8.7]{AS} again, $N_G(t)/\gen{t}\cong{}Sp_{2m-2}(q)$.  

Assume first that $P_0=1$.  If $P_1\cap{}G_0=1$, then $P_1$ is generated by 
a transvection, so $\rk([P_1,W])=2$, and $W$ does not contain a free 
$\F_2[P_2]$-module for any $P_2>P_1$.  Thus $P_1\cap{}G_0$ is a maximal 
normal unipotent subgroup of a parabolic subgroup.  So $|P_1|\ge{}q^{2m-3}$ 
(Lemma \ref{l:A5}), $|P_2|\ge{}q^{2m-2}$, and $q^{2m-2}>\rk_{\F_2}(W)=2am$ 
since we assume $(q,m)\ne(2,2)$.  Hence this case is impossible.

Next assume $P_0\ne1$ and $P_0\cap{}G_0=1$.  As noted above, $P_0$ is 
generated by a transvection (hence $\rk_{\F_q}(C_W(P_0))=2m-1$), and 
$N_G(P_0)/P_0\cong{}Sp_{2m-2}(q)$.  By the argument used in Case 1, 
$P_k/P_0\in\syl2{N_G(P_0)/P_0}$, and thus has order $q^{(m-1)^2}$.  Thus 
$|P_k/P_0|>\rk_{\F_2}(C_W(P_0))=a(2m-1)$ if $m\ge3$.  If $m=2$, then 
$N_G(P_0)/P_0\cong{}Sp_2(q)$, all nontrivial radical 2-subgroups are Sylow 
subgroups, and so there is no radical 2-chain of length 2.

Finally, assume $P_0\cap{}G_0\ne1$, and hence is a maximal normal unipotent 
subgroup of a parabolic subgroup.  Then $C_W(P_0)\le{}C_W(P_0\cap{}G_0)$ is 
a totally isotropic subspace $W_0<W$ of rank $\ell\le{}m$, and 
$N_G(P_0)/P_0\cong GL(W_0)\times N_{SO(W_0^\perp/W_0)}(P_0)/P_0$ acts on it via 
projection to the first factor.  So $P_k/P_0$ is contained in the factor 
$GL(W_0)$ since it acts faithfully on $C_W(P_0)$, $\ell\ge3$ since 
otherwise there is no radical 2-chain of length $\ge2$, and $P_k/P_0$ is a 
Sylow subgroup by the usual argument.  Thus $|P_k/P_0|=q^{\ell(\ell-1)/2}$, 
and $|P_k/P_0|>a\ell=\rk_{\F_2}(C_W(P_0))$ since $\ell\ge3$. 
\end{proof}

\begin{proof}[\textbf{Proof of Proposition \ref{ch9.1}}]  Fix a reduced 
setup $(\Gamma,S,Y)$, set $D=Z(Y)$, $V=\Omega_1(D)$, and 
$G=\Gamma/C_\Gamma(D)$, and assume $G=\gen{\calu}$ where 
	\[ \calu = \{1\ne{}P\le{}G \,|\, \textup{$P$ a quadratic best 
	offender on $D$}\}~. \]
We assume inductively that Proposition \ref{ch9.1}, and hence also 
Proposition \ref{ch8.11}, hold for all $(\Gamma^*,S^*,Y^*)$ with 
$|\Gamma^*|<|\Gamma|$.  Note that $\calu$ is a set of quadratic best 
offenders on $V$ by Lemma \ref{l:offenders}(a), and 
$O_p(G)=1$ by definition of a reduced setup.  Hence 
$G$ acts faithfully on $V$, since $C_G(V)$ is a $p$-subgroup by 
Lemma \ref{CG(V)} (so $C_G(V)\le O_p(G)$).

For each $H\le\Gamma$, let $\4H=HC_\Gamma(D)/C_\Gamma(D)$ be the image of 
$H$ in $G$. Recall that $Y=C_S(D)\in\sylp{C_\Gamma(D)}$ by definition of a 
reduced setup.  So for each $P\in\II(Y,S)$, 
$N_\Gamma(PC_\Gamma(D))=N_\Gamma(P)C_\Gamma(D)$ since 
$P\in\sylp{PC_\Gamma(D)}$. Thus
	\beqq 
	N_G(\4P)=\4{N_\Gamma(P)} \quad\textup{whenever}\quad
	Y\le P\le S .
	\label{e:4.1aa} \eeqq

Recall that $\calr=\{P\in\calf^c\,|\,J(P,D)=Y\}$.  Set 
$\calq=\II(Y,S){\sminus}\calr$.  
Thus $\calq$ is the set of all $P\in\II(Y,S)$ such that 
$\4P\cong{}P/Y$ contains at least one nontrivial best offender on $D$, 
which can be assumed quadratic by Timmesfeld's 
replacement theorem (Theorem \ref{Timmesfeld}); and $\calr$ is the set of 
all $P$ such that $\4P$ contains no such best offender.  

Set $D_0=1$.  For each $i\ge1$, set 
$D_i=\Omega_i(D)=\{g\in{}D\,|\,g^{p^i}=1\}$ and $V_i=D_i/D_{i-1}$.  Thus 
each $V_i$ is an $\F_p[G]$-module, and $(x\mapsto{}x^p)$ sends $V_i$ 
injectively to $V_{i-1}$ for each $i>0$.  

Set $k=k(p)$.  We will show that 
$\Lambda^{k}(\Out_\Gamma(R);Z(R))=0$ for each $R\in\calr$; the proposition 
then follows from Corollary \ref{c:lim*=0}.  
Here, $Z(R)=C_D(R)$ and 
$\Out_\Gamma(R)\cong{}N_\Gamma(R)/R$:  since $R\ge{}Y$ and 
$C_S(Y)=Z(Y)=D$.  So by Proposition \ref{Lambda}(d), it suffices to show,  
for each $R$ and $i$, that 
$\Lambda^{k}\bigl(N_\Gamma(R)/R;C_{D_i}(R)/C_{D_{i-1}}(R)\bigr)=0$.  Also, 
for each $i$, $C_{D_i}(R)/C_{D_{i-1}}(R)$ can be identified with an  
$N_\Gamma(R)$-invariant subgroup of $C_{V_i}(R)\le C_V(R)$. It thus 
suffices to show that 
	\beq \Lambda^{k}\bigl(N_\Gamma(R)/R;X)=0 
	\quad \textup{$\forall~ R\in\calr$,~ ~$\forall~
	N_\Gamma(R)$-invariant $X\le C_V(R)$.} 
	\eeq

Set $V_1=C_{O^p(G)}(V)$, $V_2=V_1[O^p(G),V]$, and $W=V_2/V_1$.  Thus the 
$G$-actions on $V_1$ and on $V/V_2$ factor through the $p$-group 
$G/O^p(G)$.  So by Proposition \ref{Lambda}(a,b,c), for each $R\in\calr$, 
and each $N_\Gamma(R)$-invariant $X\le{}V_1$ and $Y\le V/V_2$, 
$\Lambda^{k}(N_\Gamma(R)/R;X)=0$ and 
$\Lambda^{k}(N_\Gamma(R)/R;Y)=0$.  By the exact sequences of Proposition 
\ref{Lambda}(d), we are now reduced to showing 
	\beq \Lambda^{k}\bigl(N_\Gamma(R)/R;X)=0 
	\quad \textup{$\forall~ R\in\calr$,~ ~$\forall~
	N_\Gamma(R)$-invariant $X\le C_W(R)$.} 
	\eeq

Assume otherwise.  By Proposition \ref{Lambda-rad.ch}, 
there is a radical $p$-chain $R=P_0<P_1<P_2<\cdots<P_k$ in $\Gamma$ such 
that $C_W(P_0)$ contains a copy of the free module $\F_p[P_k/P_0]$. By 
\eqref{e:4.1aa}, $\4R=\4P_0<\4P_1<\4P_2<\cdots<\4P_k$ is a radical 
$p$-chain in $G$.  But this contradicts Proposition \ref{rad.p-ch1}, since 
by assumption, $R\in\calr$ implies that $\4R$ does not contain any 
subgroups in $\calu$. 
\end{proof}

\appendix


\newsect{Radical $p$-chains and free submodules}

We collect here some lemmas needed in the proofs in Section 4.

\begin{Lem} \label{l:A1}
Let $P$ be a $p$-group, and let $V$ be an $\F_p[P]$-module which contains a 
copy of the free module $\F_p[P]$. 
\begin{enuma} 
\item If $V=V_1\oplus\cdots\oplus V_n$ where the $V_i$ are 
$\F_p[P]$-submodules, then for some $i=1,\ldots,n$, $V_i$ contains a copy 
of $\F_p[P]$.
\item If $V$ contains an $\F_p$-basis $\bfb$ which is permuted by 
$P$, then $\bfb$ contains a free $P$-orbit.
\end{enuma}
\end{Lem}

\begin{proof} \textbf{(a) }  Set $F=\F_p[P]$ for short.  Let 
$f\:F\Right2{}V$ be a $P$-linear monomorphism.  For each $i$, let 
$\pr_i\:V\to{}V_i$ be the projection, and set $f_i=\pr_i\circ{}f$.  Choose 
$i$ such that $f_i$ sends $C_F(P)\cong\F_p$ injectively to $V_i$.  Then 
$f_i$ is injective, since otherwise $\Ker(f_i)\ne0$ would have trivial fixed 
subspace.  Thus $V_i$ contains a copy of $F$.  

\smallskip

\noindent\textbf{(b) } Write $V=V_1\oplus\cdots\oplus V_n$, where each 
$V_i$ has as basis one $P$-orbit $\bfb_i\subseteq\bfb$.  By (a), there is 
$i$ such that $V_i$ contains a copy of $\F_p[P]$, and then $\bfb_i$ is a 
free orbit.
\end{proof}

The next three lemmas involve radical $2$-subgroups and radical $2$-chains 
in symmetric and alternating groups.  

\begin{Lem} \label{l:A2}
Let $G=\Sigma_m$ or $A_m$, with its canonical action on 
$\ubm=\{1,\ldots,m\}$.  Let $P_0<P_1<\cdots<P_k$ be a radical $p$-chain in 
$G$.  Fix $0\le j\le k$, let $X_1,\ldots,X_r$ be the 
$P_j$-orbits on $\ubm$, let $\widebar{H}\le\Sigma_m$ be the subgroup of 
those permuations which send each $X_i$ to itself, and set 
$H=\widebar{H}\cap{}G$.  Thus $\widebar{H}=H_1\times\cdots\times H_r$, 
where $H_i$ is the symmetric group on $X_i$.  Let $Q_{ji}\le H_i$ be the 
image of $P_j\le H$ under the $i$-th projection, and set 
$\widebar{Q}_j=Q_{j1}\times\cdots\times Q_{jr}$ and 
$Q_j=\widebar{Q}_j\cap{}G$.  Then $Q_j=P_j$. 
\end{Lem}

\begin{proof} For each $0\le\ell\le j$, let $Q_{\ell i}\le H_i$ be the 
image of $P_\ell$ under the $i$-th projection, and set 
$\widebar{Q}_\ell=Q_{\ell1}\times\cdots\times Q_{\ell r}$ and 
$Q_\ell=\widebar{Q}_\ell\cap{}G$. We can assume by induction on $j$ that 
$Q_\ell=P_\ell$ for each $\ell<j$.  Also, $Q_{\ell i}\nsg Q_{ji}$ for each 
$i$ since $P_\ell\nsg P_j$, and so $P_\ell=Q_\ell\nsg Q_j$.

Set $G_j=N_G(P_0,\ldots,P_{j-1})$ for short. We just showed that 
$Q_j\le{}G_j$.  Each element of $\Sigma_m$ which normalizes $P_j$ must 
permute the $X_i$ and hence also normalizes $Q_j$.  Thus $N_{Q_j}(P_j)\nsg 
N_{G_j}(P_j)$, so $N_{Q_j}(P_j)=P_j$ since $P_j$ is radical in $G_j$, and 
$Q_j=P_j$ by Lemma \ref{NP>P}(a).  
\end{proof}


\begin{Lem} \label{l:A3}
Let $G=\Sigma_m$ or $A_m$, with its canonical action on 
$\ubm=\{1,\ldots,m\}$.  Let $P_0<P_1<\cdots<P_k$ be a radical $p$-chain in 
$G$ of length $k\ge2$.  Then the action of $P_k/P_0$ on $\ubm/P_0$ has no 
free orbit.
\end{Lem}

\begin{proof} We use the notation of Lemma \ref{l:A2} and its proof, with 
$j=1$.  Note that $\widebar{Q}_0$ acts trivially on $\ubm/P_0$ since each 
factor $Q_{0i}$ acts trivially on $X_i/P_0$, so 
$\widebar{Q}_1/\widebar{Q}_0$ acts on $\ubm/P_0$.  We regard $P_1/P_0$ as a 
subgroup of $\widebar{Q}_1/\widebar{Q}_0$.  Since $Q_{1i}$ acts 
transitively on $X_i$ for each $i$, $Q_{1i}/Q_{0i}$ acts transitively on 
$X_i/P_0$.

Assume there is a free $P_k/P_0$-orbit on $\ubm/P_0$, and set 
$s=|P_k/P_1|$.  After reindexing the $X_i$ if necessary, we can 
assume $P_1/P_0$ acts freely on each of its orbits $X_i/P_0$ for 
$1\le{}i\le{}s$, and that $P_k/P_1$ permutes $\{X_1,\ldots,X_s\}$ 
transitively.  Also, $P_1/P_0$ has index at most two in 
$\widebar{Q}_1/\widebar{Q}_0$, where this group splits as a product of 
subgroups $Q_{1i}/Q_{0i}$, each acting on an orbit $X_i/P_0$.  This is 
possible only if $s=|P_k/P_1|=2$, $|P_1/P_0|=2$, $|\widebar{Q}_1/\widebar{Q}_0|=4$, 
and $P_k/P_1$ exchanges $X_1$ and $X_2$ (so $|X_1|=|X_2|$).  Note that 
$\widebar{Q}_0=P_0$ since $[\widebar{Q}_1/\widebar{Q}_0:P_1/P_0]=2$.

Fix $\sigma\in{}P_1{\sminus}P_0$; thus $\sigma$ acts on $\ubm/P_0$ as a 
product of two transpositions.  For $i=1,2$, write $X_i=X'_i\amalg X''_i$, 
where $X'_i$ and $X''_i$ are $P_0$-orbits, and where $\sigma$ exchanges 
$X'_i$ with $X''_i$.  Write $\sigma=\sigma_1\cdots\sigma_r$ with 
$\sigma_i\in{}Q_{1i}$.  If $\sigma_i$ is an even permutation ($i=1,2$), 
then $(\sigma_i\times\Id_{\ubm{\sminus}X_i})\in{}P_1$, which is impossible 
by the above description of $P_1/P_0$.  Hence $\sigma_1$ and $\sigma_2$ are 
odd permutations.  

Assume $|X'_i|=\frac12|X_i|=2^t>1$ ($i=1,2$). Write 
$\sigma_i=\sigma'_i\tau_i$, where $\tau_i^2=1$, $\tau_i$ exchanges $X_i'$ 
and $X_i''$, and $\sigma'_i$ permutes $X'_i$ and is the identity on 
$X''_i$.  Then $\sigma'_i$ is odd since $\tau_i$ is even (a product of 
$2^t$ transpositions), and $\sigma_i^2=\sigma'_i\sigma''_i\in{}Q_{0i}$ 
where $\sigma''_i=\tau_i\sigma'_i\tau_i$ permutes $X''_i$ and is the 
identity on $X'_i$.  Then $\sigma'_1\sigma'_2\in{}P_0$ is the product of an 
odd permutation on $X_1$ and one on $X_2$, which is impossible since 
$\widebar{Q}_0=P_0$.  

Thus $|X'_1|=1$, so $|X_1|=|X_2|=2$, and $\sigma$ acts on the fixed subset 
$C_{\ubm}(P_0)$ as a product of two transpositions.  Hence 
$N_{\Sigma_m}(P_0,P_1)/P_1$ contains a direct factor 
$N_{\Sigma_4}(\widebar{\sigma})/\gen{\widebar{\sigma}}\cong C_2^2$ (where 
$\widebar{\sigma}=(1\,2)(3\,4)$), so $N_G(P_0,P_1)/P_1$ contains the normal 
subgroup $N_{A_4}(\widebar{\sigma})/\gen{\widebar{\sigma}}\cong C_2$.  But 
this contradicts the assumption that $P_1$ is radical in $N_G(P_0)$.  
\end{proof}


\begin{Lem} \label{l:A4}
Let $G=\Sigma_m$ or $A_m$, with its canonical action on 
$\ubm=\{1,\ldots,m\}$.  Set $V=\F_2(\ubm)$, regarded as an 
$\F_2[G]$-module.  Let $\Delta\le V$ be the $1$-dimensional submodule 
generated by the sum of the elements in $\ubm$.  Then for each radical 
$2$-subgroup $P\le{}G$, 
\begin{itemize} 
\item $C_{V/\Delta}(P)=C_V(P)/\Delta$ if $|\ubm/P|\ge3$; 
\item $\rk_{\F_2}(C_{V/\Delta}(P))\le2$ if $|\ubm/P|=2$; and 
\item $C_{V/\Delta}(P)\cong\Hom(P/\gen{\Fr(P),Q},\F_2)$ if $P$ acts 
transitively on $\ubm$ with isotropy subgroup (point stabilizer) $Q$.  
\end{itemize}
\end{Lem}

\begin{proof} The exact sequence $1\to\Delta\to V\to V/\Delta\to1$ induces 
the following exact sequence in cohomology:
	\[ 1 \Right2{} \Delta \Right3{} C_V(P) \Right3{} 
	C_{V/\Delta}(P) \Right3{} H^1(P;\Delta) 
	\Right3{\psi} H^1(P;V). \]
Here, $H^1(P;\Delta)\cong\Hom(P,\F_2)$.  If 
$x_1,\ldots,x_r\in\ubm$ are orbit representatives for the $P$-action, and 
$Q_i\le{}P$ is the stabilizer subgroup at $x_i$, then 
	\[ H^1(P;V)\cong\bigoplus_{i=1}^rH^1(P;\F_2(P/Q_i))\cong
	\bigoplus_{i=1}^r\Hom(Q_i,\F_2). \]
Also, $\psi$ is defined by restriction under these identifications, and so 
$\Ker(\psi)\cong\Hom(P/Q,\F_2)$ where $Q=\gen{\Fr(P),Q_1,\ldots,Q_r}$.  We 
thus get a short exact sequence
	\[ 1 \Right2{} C_V(P)/\Delta \Right3{} 
	C_{V/\Delta}(P) \Right3{} \Hom(P/Q,\F_2) \Right2{} 1, \]
where $\rk(C_V(P))=|\ubm/P|$.

This proves the lemma when $r=|\ubm/P|=1$.  If $r>1$, then by Lemma 
\ref{l:A2} (applied with $j=0$), $P=\widebar{P}\cap{}G$ for some 
$\widebar{P}=P_1\times\cdots\times P_r\le\Sigma_m$, where each $P_i$ 
permutes the $P$-orbit $Px_i$ transitively.  For each $i$, $P_i$ is 
contained in the stabilizer of $x_j$ for $j\ne{}i$, so $P_i\cap{}G\le{}Q$.  
Thus $[P:Q]\le2$ if $r=2$.  If $r\ge3$, then $P_iP_j\cap{}G\le Q$ for each 
$i,j\in\{1,\ldots,r\}$, and so $P=Q$ in this case.  The remaining cases of the 
lemma now follow from the above exact sequence. 
\end{proof}

We also needed lower bounds for orders of radical subgroups of 
$\Omega_{2m}^\pm(q)$.

\begin{Lem} \label{l:A5}
Let $P$ be a radical $2$-subgroup of $G=\Omega_{2m}^\pm(q)$, where $m\ge2$ 
and $q=2^a$.  Then $|P|\ge q^{2m-3}$, and $|P|\ge q^{2m-2}$ if $m\ge4$.
\end{Lem}

\begin{proof} By a theorem of Borel and Tits (see \cite[corollary]{BW} or 
\cite[13-5]{GL}), $P$ is conjugate to the maximal normal unipotent subgroup 
in a parabolic subgroup $N<G$.  Thus $P\cong O_2(N)$, $N$ contains a Borel 
subgroup and hence a Sylow 2-subgroup of $G$, and $N/O_2(N)$ is a 
semisimple group of Lie type with Dynkin diagram strictly contained in that 
of $G=D_n(q)$.  

Thus $G$ has Sylow 2-subgroups of order $q^{m(m-1)}$, while $N/O_2(N)$ must 
be contained in $D_{m-1}$ (Sylow of order $q^{(m-1)(m-2)}$), $A_{m-1}$ 
(order $q^{m(m-1)/2}$), or $A_\ell\times D_{m-\ell-1}$ (smaller order).  
Hence $|P|\ge\min\{q^{2(m-1)},q^{m(m-1)/2}\}\ge q^{2m-3}$.  When 
$G=\Omega_{2m}^-(q)$, the argument is the same, except that we only need 
consider subgroups whose Dynkin diagram is invariant under the graph 
automorphism for $D_m$.
\end{proof}


\newcommand{\bcglo}{5a}

\end{document}